\documentclass{article}

\usepackage{fullpage} 
\usepackage{parskip} 
\usepackage{tikz} 
\usepackage{amsmath}
\usepackage{amsthm}
\usepackage{hyperref}
\usepackage[normalem]{ulem} 
\usepackage{authblk} 
\usepackage{amsfonts}
\usepackage{tabularx} 
\usepackage{graphicx}
\usepackage{caption}
\usepackage{mathrsfs}
\usepackage{csquotes}
\usepackage{amssymb}
\usepackage{titling}
\usepackage{relsize}
\newcommand\numberthis{\addtocounter{equation}{1}\tag{\theequation}}

\title{\textsc{Modeling Chebyshev's Bias in the Gaussian Primes as a Random Walk}}
\author{\textsc{Daniel J. Hutama}}

\date{\textsc{July 18, 2016}}

\begin{document}

\theoremstyle{plain}
\newtheorem{thm}{Theorem}
\theoremstyle{definition}
\newtheorem{defn}{Definition}
\newtheorem*{exmp}{Example} 
\newtheorem{prpn}{Proposition}
\newtheorem{lemma}{Lemma}
\newtheorem*{corly}{Corollary}
\newtheorem*{pf}{Proof}
\newtheorem*{conj}{Conjecture}

\maketitle

\begin{abstract}
\label{abstract}
One aspect of Chebyshev's bias is the phenomenon that a prime number, $q$, modulo another prime number, $p$, experimentally seems to be slightly more likely to be a nonquadratic residue than a quadratic residue. We thought it would be interesting to model this residue bias as a ``random'' walk using Legendre symbol values as steps. Such a model would allow us to easily visualize the bias. In addition, we would be able to extend our model to other number fields. 

In this report, we first outline underlying theory and some motivations for our research. In the second section, we present our findings in the rational prime numbers. We found evidence that Chebyshev's bias, if modeled as a Legendre symbol $(\frac{q}{p})$ walk, may be somewhat reduced by only allowing $q$ to iterate over primes with nonquadratic residue (mod $4$). In the final section, we extend our Legendre symbol walks to the Gaussian primes and present our main findings. Let $\pi_1 = \alpha+\beta i$ and $\pi_2 = \beta+\alpha i$. We observed strong ($\pm$) correlations between Gaussian Legendre symbol walks for $\left[\frac{a+bi}{\pi_1}\right]$ and $\left[\frac{a+bi}{\pi_2}\right]$ where $N(\pi_1) = N(\pi_2)$ and $a+bi$ iterates over Gaussian primes in the first quadrant. We attempt an explanation of why, for some norms, the plots for $\pi_1$ and $\pi_2$ have strong positive correlation, while, for other norms, the plots have strong negative correlation.  We hope to have written in a way that makes our observations accessible to readers without prior formal training in number theory. 

\end{abstract}
\section{Introduction}
\subsection{Prime Numbers}
\label{Prime Numbers}

\begin{defn} A \textit{prime} number $p$ is any integer $p > 1$ whose divisors are only $1$ and itself. A \textit{composite} number is any integer that is not a prime number or the \textit{unit} number, $1$. \end{defn}

One of the first mathematicians to study the primes was Eratosthenes, to whom is attributed an algorithm to find all primes less than or equal to a certain value. The Sieve of Eratosthenes starts by marking all multiples of $2$ as composite, then proceeding to multiples of $3$, $5$, $7$ and so on up to $x$. 

For example, after all even numbers up to (and including) $30$ have been marked as composite, we have:

$\textbf{2}, 3, \sout{4}, 5, \sout{6}, 7, \sout{8}, 9, \sout{10}, 11, \sout{12}, 13, \sout{14}, 15, \sout{16}, 17, \sout{18}, 19, \sout{20}, 21, \sout{22}, 23, \sout{24}, 25, \sout{26}, 27, \sout{28}, 29, \sout{30}$
    
Next, we mark composite all multiples of $3$ not already marked:

$2, \textbf{3}, \sout{4}, 5, \sout{6}, 7, \sout{8}, \sout{9}, \sout{10}, 11, \sout{12}, 13, \sout{14}, \sout{15}, \sout{16}, 17, \sout{18}, 19, \sout{20}, \sout{21}, \sout{22}, 23, \sout{24}, 25, \sout{26}, \sout{27}, \sout{28}, 29, \sout{30}$
    
Next, we continue to multiples of $5$ and proceed as before, continuing until multiples of $29$:

$2, 3, \sout{4}, 5, \sout{6}, 7, \sout{8}, \sout{9}, \sout{10}, 11, \sout{12}, 13, \sout{14}, \sout{15}, \sout{16}, 17, \sout{18}, 19, \sout{20}, \sout{21}, \sout{22}, 23, \sout{24}, \sout{25}, \sout{26}, \sout{27}, \sout{28}, \textbf{29}, \sout{30}$
    
The remaining values form the set $\{2, 3, 5, 7, 11, 13, 17, 19, 23, 29\}$, which are the prime numbers less than or equal to $30$; i.e. the set of numbers less than or equal to $30$ whose divisors are only $1$ and itself.

\pagebreak

\begin{prpn}{(Fundamental Theorem of Arithmetic)} Every integer has a unique prime factorization. \end {prpn}
In other words, every integer can we expressed in a unique way as an infinite product of powers of primes:
\begin{equation}
\label{eq: fund arith}
n = 2^{\alpha_1} 3^{\alpha_2}  5^{\alpha_3} 7^{\alpha_4} \cdots = \prod p_i^{\alpha_i}
\end{equation}
where $p \in$ primes, and a finite number of $\alpha_i$ are positive integers with the rest being zero. For example, we can write $10 = 2^1 \cdot 3^0 \cdot 5^1 \cdot 7^0 \cdot 11^0 \cdots$.

\begin{prpn}{(Euclid's Theorem)} There are infinitely many prime numbers. \end{prpn}
There are many well-known proofs of Euclid's theorem. Euler's proof is as follows:\\
Let $p$ denote prime numbers and $P$ denote the set of all prime numbers. Then,
\begin{equation*}
\label{eq: eulcidsthm}
\prod_{p \in P} \sum_{\alpha \geq 0} \frac{1}{p^\alpha} = \sum_{\alpha \geq 0} \frac{1}{2^\alpha} \cdot \sum_{\alpha \geq 0} \frac{1}{3^\alpha} \cdot \sum_{\alpha \geq 0} \frac{1}{5^\alpha} \cdot \sum_{\alpha \geq 0} \frac{1}{7^\alpha} \cdots = \sum_{\alpha_1 , \alpha_2, \alpha_3, \ldots \geq 0} \frac{1}{2^{\alpha_1} 3^{\alpha_2} 5^{\alpha_3} \cdots}
\end{equation*}

However, by \eqref{eq: fund arith}, we know that every integer can be written uniquely as a product of primes. Thus, we can rewrite our equation as:
\begin{equation}
\label{eq: eulersproof}
\prod_{p \in P} \sum_{\alpha \geq 0} \frac{1}{p^\alpha} = \sum_{\alpha_1 , \alpha_2, \alpha_3, \ldots \geq 0} \frac{1}{2^{\alpha_1} 3^{\alpha_2} 5^{\alpha_3} \cdots} = \sum_{n} \frac{1}{n}
\end{equation}

We then recognize the right hand side of \eqref{eq: eulersproof} as the harmonic series. Because of the divergence of the harmonic series, we know our product must be infinite as well. Since each term of our product is a finite number, there must be an infinite number of terms for the product to be infinite.

Euler also proved a stronger version of the divergence of the harmonic series, in which he shows the sum of reciprocals of primes also diverges \cite{prime reciprocal divergence}. We will use this fact in a later proof.
\begin{equation}
\label{eq:sum of reciprocals}
\sum_{p \in P} \frac{1}{p} = \infty
\end{equation}

\subsection{Arithmetic Progressions}
\label{Arithmetic Progressions}

The Sieve of Eratosthenes is effective because of the simplicity of identifying multiples of a number. For example, it is easy to identify all numbers of the form $3n$ (which is the set $\{3, 6, 9, 12, 15, 18, \ldots \}$ for $n\geq 1$) as multiples of $3$, and subsequently mark them as composite (with the exception of the first element). However, what happens if we change the starting value of the set, while keeping the distance between elements the same?

\begin{defn} We call a sequence of numbers with constant difference between terms an \textit{arithmetic progression}.\end{defn}

For example, consider all numbers of the form $3n+2$ and $3n+1$, which represent the sets $\{2, 5, 8, 11, 14, 17 \ldots \}$ and $\{1, 4, 7, 10, 13, 16, \ldots \}$ respectively. Both sets of numbers are arithmetic progressions with a difference of $3$. 

The reader might then inquire: 
\begin{itemize}
\item Between $3n+2$ and $3n+1$, which arithmetic progression contains more primes up to a value $x$? In other words, if we consider the count of primes in each progression as a race, which team is in the lead at a given $x$? 
\item Can we extend Euclid's Theorem to primes in arithmetic progressions? In other words, do arithmetic progressions contain infinitely many primes?
\item What is the distribution of primes in these progressions?
\end{itemize}
To answer these questions, we must first introduce a few tools to give our analysis some sophistication.

\subsection{Euclidean Algorithm, Euler's Totient Function, and Modulo}
\label{Euclidean Algorithm}

\begin{defn} An integer $a \not = 0$ \textit{divides} another integer $b$ if there exists another integer $c$, such that $b=ac$. We denote that $a$ divides $b$ with $a \vert b$. \end{defn} 

\begin{defn}
Pick two integers $a$ and $b$. An integer $c$ such that  $c \vert a$ and $c \vert b$ is said to be a \textit{common divisor} of $a$ and $b$. If there exists another integer $d \geq c$ that also divides $a$ and $b$, we say that $d$ is the \textit{greatest common divisor} of $a$ and $b$. We denote this by $\gcd(a,b)=d$. 
\end{defn}

\begin{prpn}
Let $a$ and $b$ be integers. The Euclidean Algorithm allows us to compute the greatest common divisor of $a$ and $b$; i.e. it allows us to find the largest number that divides both $a$ and $b$, leaving no remainder. The algorithm is as follows:
\end{prpn}
\begin{align*}
a & =bq_0+r_0 && \text{for} && 0<r_0<b \\
b & =r_0 q_1+r_1 && \text{for} && 0<r_1<r_0 \\ 
r & =r_1 q_2+r_2 && \text{for} && 0<r_2<r_1   \\ 
& && \ldots &&& \\
r_{k-1} & = r_k q_{k+1}+r_{k+1} && \text{for} && 0<r_{k+1}<r_k \\
r_k & =r_{k+1} q_{k+2}+ 0
\end{align*}

Then $\gcd(a, b) = r_{k+1}$. For example, to find $\gcd(6188,4709)$, we apply the Euclidean Algorithm as follows:
\begin{align*}
6188 & = 4709 \cdot 1 + 1479 \\
4709 & = 1479\cdot 3 + 272 \\
1479 & = 272 \cdot 5 + 119 \\
272 & = 119 \cdot 2 + 34 \\
119 & = 34 \cdot 3 + 17 \\
34 & = 17 \cdot 2 \\
17 & = \gcd(6188,4709) 
\end{align*}

\begin{defn}
$a$ and $b$ are said to be \textit{relatively prime}, or \textit{coprime} if $\gcd(a,b) =1$. 
\end{defn}

Two prime numbers, $p$ and $q$, will always be coprime to each other. A composite number, $a$, will be coprime to prime number, $p$, if and only if $a$ is not a multiple of $p$.

\begin{defn}
\textit{Euler's totient function}, denoted $\phi(n)$, counts the number of \textit{totatives} of $n$, i.e. the number of (positive) integers up to $n$ that are coprime to $n$. 
\end{defn}

For example, $\phi(10) = \#\{1,3,7,9\} = 4$. In this example, the numbers $1$, $3$, $7$, and $9$, are the totatives of $10$. For a prime number $p$, $\phi(p) = \#\{1,2,\ldots,p-1\} = p-1$ since all integers $< p$ are also coprime to $p$.

\begin{defn}
We say that $a$ is \textit{congruent} to $r$ \textit{modulo} $b$ if $b \vert a-r$. We write this relation as $a \equiv r \pmod{b}$ 
\end{defn}
In other words, we say that $a \equiv r$ (mod $b$) if $r$ is the remainder when $a$ is divided by $b$. For example, when $9$ is divided by $7$, the remainder is $2$. In other words, $9 \equiv 2$ (mod $7$). This concept allows us to conveniently refer to arithmetic progressions by their congruences modulo $a$. For instance, we can refer to the progression $4n+3$ as the set of all integers congruent to $3$ (mod $4$). Furthermore, we can refer to all primes in the progression $4n+3$ as the set of primes congruent to $3$ (mod $4$).

\begin{corly}
Let $\mathbb{Z}$ denote the set of all integers. The modulo operation allows us to define a quotient ring, $\mathbb{Z}/n\mathbb{Z}$, which is the ring of integers modulo $n$. \end{corly}
For example, the set of all integers modulo $6$ repeats as $\{\ldots, 1,2,3,4,5,0,1,2,3,4,5,\ldots\}$. The unique elements of this set are $\{0,1,2,3,4,5\}$, which is the ring $\mathbb{Z}/6\mathbb{Z}$. We say that an element $u$ in $\mathbb{Z}/n\mathbb{Z}$ is a \textit{unit} in the ring if there exists a multiplicative element $v$, such that $uv = vu = 1$. We denote the group of units as $(\mathbb{Z}/n\mathbb{Z})^{\times}$. 

The group $(\mathbb{Z}/n\mathbb{Z})^{\times}$ has $\phi(n)$ elements, which are the totatives of $n$. For example, for the ring $\mathbb{Z}/6\mathbb{Z}$, the group of units, $(\mathbb{Z}/6\mathbb{Z})^{\times}$ is given by the totatives of 6: \{1,5\}. We notice that $1$ and $5$ are both units in $\mathbb{Z}/6\mathbb{Z}$ since $1\equiv 1 \pmod 6$ and $5\cdot 5 \equiv 1 \pmod 6$.
Thus for a prime number $p$, the group $(\mathbb{Z}/p\mathbb{Z})^{\times}$ has $p-1 = \phi(p)$ elements. 

\subsection{The Prime Number Theorem and Dirichlet's Theorem on Arithmetic Progressions}
\label{Prime Number Theorem}

Let $\pi(x)$ denote the number of primes up to $x$. \\
\begin{prpn}
Gauss's Prime Number Theorem (PNT), which Hadamard and Vall\`{e}e-Poussin proved independently in 1896, states that $\pi(x)$ behaves asymptotically to $x/\log(x)$\footnote{ $\log(x)$ here is actually the natural log of $x$, but we wish to use the same notation as in our references}
\end{prpn}
Put another way:
\begin{equation}
\label{eq:PNT}
\lim_{x \to \infty} \frac{\pi(x)}{x/\log(x)} = 1
\end{equation}

Thus for an arbitrarily large value of $x$, one can expect $\pi(x)$ to be close to $x/\log(x)$, with some error term. One might next wonder about approximating the count of primes within an arithmetic progression. One way of intuitively approaching this problem is by viewing the set of all positive integers as a union of arithmetic progressions. For example, if we consider the arithmetic progressions with a difference of $3$ between elements in each set, we have the three progressions: 
\begin{align*}
& \{3n+1 &&\text{for} &&n \in \mathbb{N}_0\} &=&&\{1,4,7,10,13,16,\ldots\} \\
& \{3n+2 &&\text{for} &&n \in \mathbb{N}_0\} &=&&\{2,5,8,11,14,17,\ldots\} \\
& \{3n &&\text{for} &&n \in \mathbb{N}_1\} &=&&\{3,6,9,12,15,18,\ldots\}
\end{align*}
Combining these three sets will yield the set of all positive integers. Since each element in the third set is a multiple of $3$, and thus a composite number, we can ignore this set and only consider the first two. We can then expect the primes to be split approximately equally between $3n+1$ and $3n+2$. Similarly, for a difference of $4$ between elements in each set, primes would be split approximately evenly between $4n+1$ and $4n+3$. 

Thus applying our intuition to \eqref{eq:PNT}, we arrive at:
\begin{thm} 
\label{dirichletsthm}
(Dirichlet's Theorem on Arithmetic Progressions)
If $\gcd(a,b) = 1$, there are infinitely many primes congruent to $b$ modulo $a$.  In addition, for progressions of the form $an+b$, the primes will be split among $\phi(a)$ different progressions. In other words, the proportion of primes in a progression with increment $a$ is $\frac{1}{\phi(a)}$. 
\end{thm}
\begin{equation}
\label{eq:PNTforAPs}
\lim_{x \to \infty} \frac{\pi(x;a,b)}{x/(\phi(a) \cdot \log(x))} = 1
\end{equation}
For example, the progression $5n+1$ holds one-fourth of primes ($\phi(5)=4$), and we write:
\begin{equation*}
\lim_{x \to \infty} \frac{\pi(x;5,1)}{x/(\phi(5) \cdot \log(x))} = 1
\end{equation*}

\pagebreak

The complete proof of Dirichlet's Theorem is quite lengthy, but excellently shown by Pete L. Clark \cite{dirichlets_theorem} and Austin Tran \cite{dirichlets theorem orthogonality}. Here, we only briefly introduce important concepts from analytic number theory and highlight crucial points of the proof as shown by Clark and Tran. For readers not familiar with analytic number theory, this section may be particularly difficult. Nevertheless, we encourage the reader on.

\begin{defn}
A \textit{Dirichlet Character} modulo $a$ is a function $\chi$ on the units of $\mathbb{Z}/a\mathbb{Z}$ that has the following properties:
\end{defn}

\begin{itemize}
\item $\chi$ is periodic modulo $a$, i.e. $\chi(b) = \chi(b+a)$ for $b \in \mathbb{N}$.
\item $\chi$ is multiplicative, i.e. $\chi(b) \cdot \chi(c) = \chi(bc)$.
\item $\chi(1) = 1$.
\item $\chi(b) \not = 0$ if and only if $\gcd(a, b) = 1$.
\end{itemize}

We say that a character is \textit{principal} if its value is $1$ for all arguments coprime to its modulus, and $0$ otherwise. We denote the principal character modulo $a$ as $\chi_0$. Note that the principal character still depends on $a$.

\begin{exmp}
Consider the Dirichlet characters modulo 3. We have $\chi(1) = 1$ and 
$\chi(3) = 0$ by properties stated above.  Using the multiplicativity and periodicity of $\chi$ we note that $(\chi(2))^2 = \chi(2)\cdot \chi(2) = \chi(1) = 1$. This implies that $\sqrt{(\chi(2))^2} = \chi(2) = \pm 1$. If $\chi(2) = 1$, then $\chi = \chi_0$ is a principal character by definition. On the other hand, we use $\chi_1$ to denote the character for when $\chi(2) = -1$. We note that $\chi_1$ also satisfies all necessary properties to be a Dirichlet character, but is not a principal character. 
\end{exmp}

\begin{prpn}
Let $X(a)$ denote the set of all Dirichlet Characters modulo $a$. $X(a)$ is a group with multiplication and an identity element given by the principal character $\chi_0$ modulo $a$. In addition, the following orthogonality relation holds (orthogonality of characters):
\end{prpn}
\begin{equation*}
\sum_{\chi \pmod a} =  \begin{cases}
        1  & \text{if $b \equiv 1$ (mod $a$)}, \\
        0 & \text{otherwise}
        \end{cases}
\end{equation*}
(A proof of the orthogonality of characters is nicely shown by A. Tran in \cite{dirichlets theorem orthogonality}). 

\begin{corly}
The values of a character $\chi$ are either $0$ or the $\phi(a)^{\text{th}}$ roots of unity.
\end{corly}

Recall that if $\chi(b) \not = 0$, then $\gcd(a, b) = 1$. If order of the group is $\phi(a)$, then $\chi(b)^{\phi(a)}$ is principal, so $\chi(b)^{\phi(a)} = 1$.  Thus, $\chi(b) = e^{\frac{2\pi i \nu}{\phi(a)}}$ for $\nu \in \mathbb{N}$.

\begin{defn}
A \textit{Dirichlet L-series} is a function of the form:
\end{defn}
\begin{equation*}
L(\chi, s) = \sum_{n=1}^{\infty}\frac{\chi(k)}{n^s}
\end{equation*}
where $s$ is a complex variable with Re($s$)$>1$. 
\begin{prpn}
\label{prpn: euler product}
The Dirichlet L-function can be also expressed as an Euler product as follows (A proof can be found in \cite{euler product reference}):
\end{prpn}
\begin{equation}
\label{eq: L series Euler rep}
L(\chi, s) = \prod_p \left(1 - \frac{\chi(p)}{p^s}\right)^{-1}
\end{equation}

We introduce an intermediate theorem necessary for the  proof of theorem \ref{dirichletsthm}:
\begin{thm}
\label{prop: nonvanishing}
\textit{Dirichlet's Non-vanishing Theorem} states that $L(\chi,1) \not = 0$ if $\chi$ is not a principal character. 
\end{thm}

Here, we will only highlight crucial sections of the proof of Dirichlet's non-vanishing theorem (as shown by J.P. Serre). A more complete proof of Theorem \ref{prop: nonvanishing} can be found in \cite{jp serre}.

Let $a$ be a fixed integer $\geq 1$. If $p \nmid m$, we denote the image of $p$ in $(\mathbb{Z}/a\mathbb{Z})^\times$ by $\overline{p}$. In addition, we use $f(p)$ to denote the order of $p$ in $(\mathbb{Z}/a\mathbb{Z})^\times$; i.e. $f(p)$ is the smallest integer $f$ such that $p^f \equiv 1 \pmod a$. We let $g(p) = \frac{\phi(a)}{f(p)}$. This is the order of the quotient of $(\mathbb{Z}/a\mathbb{Z})^\times$ by the subgroup $(\overline{p})$ generated by $p$. 

\begin{lemma}
\label{lemma1}
For $p \nmid a$, we have the identity:
\end{lemma}
\begin{equation*}
\prod_{\chi \in X(a)}(1 - \chi(p)T) = (1-T^{f(p)})^{g(p)}
\end{equation*}
For the derivation of lemma \ref{lemma1}, we let $\mu_{f(p)}$ denote the set of $f(p)^{th}$ roots of unity. We then have the identity:
\begin{equation}
\label{eq:nonvanish identity}
\prod_{w \in \mu_{f(p)}}(1-wT) = 1-T^{f(p)}
\end{equation}
For all $w\in\mu_{f(p)}$, there exists $g(p)$ characters $\chi \in X(a)$ such that $\chi(\overline{p}) = w$. This fact, together with \eqref{eq:nonvanish identity}, brings us to lemma \ref{lemma1}. 

We now define a function $\zeta_a(s)$ as follows:
\begin{equation*}
\zeta_a(s) = \prod_{\chi \in X(a)}L(\chi, s)
\end{equation*}
We continue by replacing each $L(\chi,s)$ in the product by its product expansion as in \eqref{eq: L series Euler rep}, and then applying lemma \ref{lemma1} with $T=p^{-s}$.
\begin{prpn} 
\label{prpn: zeta expansion}
We can then represent the product expansion of $\zeta_a(s)$ as follows:
\end{prpn}
\begin{equation*}
\zeta_a(s) = \prod_{p \nmid a}\dfrac{1}{\left(1-\dfrac{1}{p^{f(p)s}}\right)^{g(p)}}
\end{equation*}
We note that this is a Dirichlet series with positive integral coefficients converging in the half plane $Re(s) > 1$.

We now wish to show (a) that $\zeta_a(s)$ has a simple pole at $s=1$ and (b) that $L(\chi, 1) \not = 0$ for all $\chi \not = \chi_0$. The fact that $L(1, s)$ has a simple pole at $s=1$ implies the same for $\zeta_a(s)$. Thus, showing (b) would imply (a).

Suppose for contradiction that $L(\chi,1) = 0$ for $\chi \not = \chi_0$. Then $\zeta_a(s)$ would be holomorphic at $s=1$, and also for all $s$ with $Re(s)>0$. Since by proposition \ref{prpn: zeta expansion}, $\zeta_a(s)$ is a Dirichlet series with positive coefficients, the series would converge for all $s$ in that domain. However, this cannot be true. We show this by expanding the $p^{th}$ factor of $\zeta_a(s)$ as follows:
\begin{equation*}
\dfrac{1}{(1-p^{-f(p)s})^{g(p)}} = (1+p^{-f(p)s} + p^{-2f(p)s} + p^{-3f(p)s}+\ldots)
\end{equation*}
We then ignore crossterms with negative contribution to arrive at an upper bound:
\begin{equation*}
1 + \dfrac{1}{p^{\phi(a)s}}+\dfrac{1}{p^{2\phi(a)s}}+\dfrac{1}{p^{3\phi(a)s}} + \ldots
\end{equation*}

Multiplying over $p$, it follows that $\zeta_a(s)$ has all its coefficients greater than the series:
\begin{equation}
\label{eq:divergence at phi a}
\sum_{n | \gcd(a,n)=1}\dfrac{1}{n^{\phi(a)s}} 
\end{equation}
Evaluating equation \eqref{eq:divergence at phi a} at $s=\frac{1}{\phi(a)}$, we finish the proof of theorem \ref{prop: nonvanishing} by arriving at the following divergent series:
\begin{equation*}
\sum_{n | \gcd(a,n)=1}\dfrac{1}{n}.
\end{equation*}

We now proceed with the proof of Dirichlet's Theorem.

\textit{Proof of Theorem 1}. Let $X(a)$ denote the group of Dirichlet characters modulo $a$. We then fix $\gcd(a,b)=1$ as stated in Dirichlet's Theorem. In addition, we let $\Psi$ denote the set of prime numbers $p \equiv b \pmod a$. Our goal is to show that $\Psi$ is an infinite set.

We wish to consider a function similar to the one in \eqref{eq: eulersproof}.  We define:
\begin{equation}
\label{eq:congruent primes recip}
P_b(s) := \sum_{p\in \Psi} \frac{1}{p^s}
\end{equation}
In particular, we wish to show that the function $P_b(s)$ approaches $\infty$ as $s$ approaches $1$. This would imply infinitely many elements in $\Psi$. We also define $\theta_b$ to be the characteristic function of the congruence class $b \pmod a$. In other words: 
\begin{equation*}
\theta_b(n) =  \begin{cases}
        1  & \text{if $n \equiv b$ (mod $a$)}, \\
        0 & \text{otherwise}
        \end{cases}
\end{equation*}

Note that $\theta_b$ is periodic modulo $a$ and is $0$ when $\gcd(n,a)>1$. 

Using this characteristic function, we wish to express $P_b(s)$ as a sum over all primes:
\begin{equation*}
P_b(s) = \sum_{p\in P} \frac{\theta_b(p)}{p^s}
\end{equation*}

\begin{lemma}
\label{lemma: theta}
For all $n \in \mathbb{Z}$, we have:
\end{lemma}
$$\theta_b = \sum_{\chi \in X(a)}\frac{\chi(b^{-1})}{\phi(a)}\chi(n)$$

\textit{Proof of Lemma \ref{lemma: theta}}. Using the multiplicative property of the Dirichlet character:
$$\theta_b = \frac{1}{\phi(a)}\left(\sum_{\chi \in X(a)}\chi(b^{-1}n)\right)$$
By our orthogonality relation, the summation term becomes $\phi(a)$ if $b^{-1}n = 1$ (i.e. if $n\equiv b \pmod a$) and zero otherwise. The result is exactly $\theta_b$.

Applying Lemma \ref{lemma: theta} to \eqref{eq:congruent primes recip}, we arrive at:
\begin{equation}
\label{eq:dirichletproof2}
P_b(s) = \sum_{\chi \in X(a)}\frac{\chi(b^{-1})}{\phi(a)}\sum_p\frac{\chi(p)}{p^s}
\end{equation}

We recognize the second summation term as reminiscent of the Dirichlet series we defined earlier. We will come back to this equation later. 

Consider the convergent Taylor series expansion of $\log(1-z)$ for $\vert z \vert < 1$
\begin{equation}
\label{eq:complex log}
\log(1-z) = - \sum_{n=1}^{\infty}\frac{z^n}{n}
\end{equation}
 In addition, consider the Euler product representation of our Dirichlet series in \eqref{eq: L series Euler rep}. Applying logarithms, we get:
 \begin{equation}
 \label{eq: log of euler prod of L series}
 \log(L(\chi, s)) = - \sum_p \log\left(1-\frac{\chi(p)}{p^s}\right)
\end{equation}

Combining \eqref{eq:complex log} and \eqref{eq: log of euler prod of L series}, we have:
\begin{equation}
\label{eq: log euler with complex log}
  \log(L(\chi, s)) = \sum_p\sum_n \frac{1}{n}\left(\frac{\chi(p)}{p^s}\right)^n
\end{equation}

The right side of \eqref{eq: log euler with complex log} is absolutely convergent for Re($s$) $> 1$, and is therefore an analytic function on that half plane. We now denote the right hand side of \eqref{eq: log euler with complex log} as $\mathit{l}(\chi, s)$.

\begin{lemma}
In the half plane with Re($s$)$>1$, $e^{\mathit{l}(\chi,s)} = L(\chi, s)$.
\end{lemma}
The proof of Lemma 2 is shown in \cite{dirichlets theorem orthogonality}.

We now split $\mathit{l}(\chi, s)$ into two parts. The first part will be for the sums when $n=1$, and the second part will be for the sums when $n>1$. We denote these as $\mathit{I}(\chi, s)$ and $\mathit{R}(\chi, s)$ respectively. Symbolically, 
$$ \mathit{l}(\chi, s) = \mathit{I}(\chi,s) + \mathit{R}(\chi, s)$$ $$\mathit{I}(\chi,s) = \sum_p \frac{\chi(p)}{p^s}, \mathit{R}(s,\chi) = \sum_{n \geq 2}\sum_{p}\frac{\chi(p)^n}{np^{ns}}$$

We now note that we can write $P_b(s)$ from \eqref{eq:dirichletproof2} as:

\begin{equation}
\label{eq: split equation applied to l series log}
P_b(s) = \sum_{\chi \in X(a)} \frac{\chi(b^{-1})}{\phi(a)} \mathit{I}(\chi, s)
\end{equation}

\begin{lemma} $\mathit{R}(\chi, s)$ is bounded when $s=1$ (Recall, that we wish to show that $P_b(s) \rightarrow \infty$ as $s \rightarrow 1$). \end{lemma}
 This can be shown by comparing $\mathit{R}(\chi, s)$ to the well-known Basel problem:

\begin{equation*}
\vert \mathit{R}(\chi, 1) \vert \leq \sum_{n\leq 2} \sum_{p} \frac{1}{np^n} \leq \sum_{p} \sum_{n \leq 2} \frac{1}{p^n} \leq 2  \sum_n \frac{1}{n^2} = \frac{2\pi^2}{6}
\end{equation*}

Since we know that $\mathit{R}(\chi, 1)$ is bounded, we can ignore it as it will not help us in showing that $P_b(s)$ diverges as $s \rightarrow 1$.

We now wish to split our summation from  \eqref{eq: split equation applied to l series log} into an expression with only principal characters, and a sum over non-principal characters. Recall that 
a principal character $\chi_0(n) = 1$ for $\gcd(n,a) =1$, and 0 otherwise. 

\begin{equation*}
P_b(s) = \sum_{\chi \in X(a)}\frac{\chi(b^{-1})}{\phi(a)}\mathit{I} (\chi, s)
\end{equation*}
\begin{equation*}
= \frac{\chi_0(b^{-1})}{\phi(a)}\mathit{I}(\chi_0, s) +  \sum_{\chi \not = \chi_0}\frac{\chi(b^{-1})}{\phi(a)}\mathit{I} (\chi, s)
\end{equation*}
\begin{equation}
\label{eq: almost done with proof}
P_b(s) = \frac{1}{\phi(a)}\sum_{p \nmid a} \frac{1}{p^s} + \sum_{\chi \not = \chi_0}\mathit{l}(\chi, s)
\end{equation}

We know that $a$ will have a finite number of prime divisors. This fact, together with equation \eqref{eq:sum of reciprocals}, tells us that the first term in \eqref{eq: almost done with proof} is unbounded. All that remains is to show that the second summation in \eqref{eq: almost done with proof} is bounded as $s \rightarrow 1$. Doing so will show that the primes (mod $a$) will fall into one of the $\phi(a)$ congruence classes as claimed in theorem \ref{dirichletsthm}. To do this, we must use Dirichlet's non-vanishing theorem (theorem \ref{prop: nonvanishing}). Recall that $L(\chi, 1) \not = 0$ if $\chi$ is not a principal character. Thus:
\begin{equation*}
\label{eq: almost last eq in proof of dirichlets theorem}
L(\chi, s) = \lim_{s\rightarrow 1}L(\chi, s) = \lim_{s \rightarrow 1}e^{\mathit{l}(\chi, s)}
\end{equation*}

Since logarithms of an analytic function differ only by multiples of $2\pi i$, $\mathit{l}(\chi, s) = \log L(\chi, s)$ always remains bounded as $s \rightarrow 1$. As a result, the contribution to $P_b(s)$ from non-principal Dirichlet characters remains bounded, while the contribution from principal characters is unbounded. $P_b(s)$ itself is then unbounded as $s \rightarrow 1$. In conclusion, we have:
\begin{equation*}
\sum_{p \in \Psi}\frac{1}{p^s} = \lim_{s \rightarrow 1}P_b(s) = \infty
\end{equation*}

Thus, there must be infinitely many elements in $\Psi$, i.e. there are infinitely many primes congruent to $b$ modulo $a$ for $\gcd(a, b) = 1$.

\subsection{Chebyshev's Bias, Quadratic Residue, and the Legendre Symbol}
\label{Chebyshev's Bias}

As quite thoroughly shown by A. Granville and G. Martin in their paper, Prime Number Races \cite{prime_number_races}, when we ``race'' progressions, some progressions hold the lead for an overwhelming majority of the time. For example, in the mod $4$ race of $4n+1$ against $4n+3$, the bias is as much as $99.59\%$ in favor of the $4n+3$ team!

This bias, first observed by Chebyshev in 1853, is attributed to primes in the $4n+1$ progression being \textit{quadratic residues} modulo $4$. As noted by Terry Tao \cite{consecutive_biases}:
\begin{displayquote}
...Chebyshev bias asserts, roughly speaking, that a randomly selected prime $p$ of a large magnitude $x$ will typically (though not always) be slightly more likely to be a quadratic non-residue modulo $q$ than a quadratic residue, but the bias is small (the difference in probabilities is only about $\mathit{O}(\frac{1}{\sqrt{x}})$ for typical choices of $x$)
\end{displayquote}

\begin{defn}
Let $p$ be an odd prime number\footnote{Restriction is such that the Legendre symbol will be defined for any $p$.} . We say that a number $a$ is a \textit{quadratic residue} (QR) modulo $p$ if there exists an element $x$ in the set of totatives of $p$, such that $x^2 \equiv a$ (mod $p$).
\end{defn}

(Note: $p$ does not necessarily need to be prime for the definition of quadratic residues. However, as we will see later, the modulus must be prime for our Legendre symbol model to work. Thus, we restrict our study to only prime moduli).

For example, let us consider the set of totatives of $7$, which is the set \{1,2,3,4,5,6\}:
\begin{table}[h]
\centering
\caption{Quadratic Residues (mod 7)}
\label{QRtable}
\begin{tabular}{|c|c|c|c|l}
\cline{1-4}
$x$ & $x^2$ & $x^2$ (mod 7) & Conclusion &  \\ \cline{1-4}
1 & 1                   & 1                           & 1 is a QR (mod 7) &  \\ \cline{1-4}
2 & 4                   & 4                           & 4 is a QR (mod 7) &  \\ \cline{1-4}
3 & 9                   & 2                           & 2 is a QR (mod 7) &  \\ \cline{1-4}
4 & 16                  & 2                           & 2 is a QR (mod 7) &  \\ \cline{1-4}
5 & 25                  & 4                           & 4 is a QR (mod 7) &  \\ \cline{1-4}
6 & 36                  & 1                           & 1 is a QR (mod 7) &  \\ \cline{1-4}
\end{tabular}
\end{table}

In this example, $1$ is a quadratic residue since both $1^2$ and $6^2$ are congruent to $1$ (mod 7). In addition, $4$ is a quadratic residue since $2^2$ and $5^2$ are congruent to $4$ (mod 7), and $2$ is a quadratic residue since $3^2$ and $4^2$ are congruent to $2$ (mod 7). Note the symmetry of quadratic residues when ordered by $x$.

We now might like a convenient notation to quantify the notion of quadratic residues.

 \pagebreak
 
\begin{defn}
The Legendre symbol separates an integer $a$ into three classes, depending on its residue modulo an odd prime $p$. 
\end{defn}

\begin{equation*}
\left(\frac{a}{p}\right) =  \begin{cases}
        1  & \text{if $a$ is a quadratic residue (mod $p$)}, \\
        -1 & \text{if $a$ is a nonquadratic residue (mod $p$)}, \\
        0 & \text{if $a \equiv 0$ (mod $p$)}
        \end{cases} 
\end{equation*} 

Note: the Legendre symbol is only defined for $p$ being an odd prime number. If $a$ is a prime number $\not = p$, the Legendre symbol will never be $0$ (since two different prime numbers will be coprime) We know by Theorem \ref{dirichletsthm} that the residues of $a \pmod p$ are then equally distributed among congruence classes in $\{1,2,3 \ldots, p-1\}$.

Continuing with our definition, we introduce several properties of the Legendre symbol:

\begin{itemize}
\item The Legendre symbol is periodic on its top argument modulo $p$. In other words, if $a \equiv b \pmod p$, then
\begin{equation*}
\left(\frac{a}{p}\right) = \left(\frac{b}{p}\right)
\end{equation*}
\item The Legendre symbol is multiplicative on its top argument, i.e.
\begin{equation*}
\left(\frac{a}{p}\right)\left(\frac{b}{p}\right) = \left(\frac{ab}{p}\right) 
\end{equation*}
\item The product of two squares is a square. The product of two nonsquares is a square. The product of a square and a nonsquare is a nonsquare. This can be expressed as follows:
\begin{align*}
&\text{Two squares}:&1 \cdot 1 &= 1 \\
&\text{Two nonsquares}: &-1 \cdot -1 &= 1 \\
&\text{Square and nonsquare}:& 1 \cdot -1 &= -1
\end{align*}
\item The Legendre symbol can also be defined equivalently using Euler's criterion as:
$$\left(\frac{a}{p}\right) \equiv \mathlarger{a}^{(p-1)/2} \pmod p$$
\end{itemize}

\begin{prpn} 
\label{prpn: quadratic reciprocity}
(Law of Quadratic Reciprocity) For $p$ and $q$ odd prime numbers:
\end{prpn}
\begin{equation*}
\left(\frac{q}{p}\right) = \mathlarger{(-1)}^{\frac{p-1}{2}\frac{q-1}{2}} \left(\frac{p}{q}\right)
\end{equation*}
The Law of Quadratic Reciprocity \cite{QR stein} has several supplements for different values of $a$. Here, we only introduce the first two supplements without proof. For $x$ in the set of totatives of $p$:

\begin{enumerate}
\item $x^2 \equiv -1 \pmod p$ is solvable if and only if $p \equiv 1 \pmod 4$.
\item $x^2 \equiv 2 \pmod p$ is solvable if and only if $p \equiv \pm 1 \pmod 8$.
\end{enumerate}

These supplements can be expressed equivalently as follows:
\begin{enumerate}
\item \begin{equation*}
\left(\frac{-1}{p}\right) = \mathlarger{(-1)}^{\frac{p-1}{2}} = 
		\begin{cases}
        1  & \text{if $p \equiv 1$ (mod $4$)}, \\
        -1 & \text{if $p \equiv 3$ (mod $4$)} 
        \end{cases} 
\end{equation*}
\item \begin{equation*}
\left(\frac{2}{p}\right) = \mathlarger{(-1)}^{\frac{p^2-1}{8}} =
		\begin{cases}
        1  & \text{if $p \equiv 1$, $p \equiv 7$ (mod $8$)}, \\
        -1 & \text{if $p \equiv 3$, $p \equiv 5$ (mod $8$)} 
        \end{cases} 
\end{equation*}
\end{enumerate}

\pagebreak
Continuing with our example for $a$ in $(\mathbb{Z}/7\mathbb{Z})^{\times}$, we have:
\begin{table}[h]
\centering
\caption{Legendre Symbols (mod 7)}
\label{LStable}
\begin{tabular}{|c|c|c|c|c|c|c|}
\cline{1-7}
$a$                        & $1$ & $2$ & $3$  & $4$ & $5$  & $6$          \\ \cline{1-7}
$\left(\frac{a}{7}\right)$ & $1$ & $1$ & $-1$ & $1$ & $-1$ & $-1$        \\ \cline{1-7}
\end{tabular}
\end{table}

\begin{prpn}
In general, Chebyshev's bias suggests that, in a race between $\alpha n + \beta_1$ and $\alpha n+\beta_2$, the progression in which $\beta_i$ is a nonquadratic residue (mod $\alpha$) will likely contain more primes up to $x$. 
\end{prpn}

For instance, when racing $1$ (mod $3$) against $2$ (mod $3$), we observe that $2$ (mod $3$) almost always has more primes up to $x$. Indeed, $1$ is a quadratic residue (mod $3$), and $2$ is a nonquadratic residue (mod $3$). 

\begin{table}[h]
\centering
\caption{Count of Primes in the mod 3 Race}
\label{mod3race}
\begin{tabular}{|c|c|c|l}
\cline{1-3}
$x$ & Primes in $3n+1$ up to $x$  & Primes in $3n+2$ up to $x$  \\ \cline{1-3}
$10^1$ & 1                    & 2                            &    \\ \cline{1-3}
$10^2$ & 11                   & 13                           &    \\ \cline{1-3}
$10^3$ & 80                   & 87                           &    \\ \cline{1-3}
$10^4$ & 611                  & 617                          &    \\ \cline{1-3}
$10^5$ & 4784                 & 4807                         &    \\ \cline{1-3}
$10^6$ & 39231                & 39266                        &    \\ \cline{1-3}
\end{tabular}
\end{table}

Despite the apparent domination by the $2$ (mod $3$) team, a theorem from J.E. Littlewood (1914) asserts that there are infinitely many values of $x$ for which the $1$ (mod $3$) team is in the lead (of course, this theorem applies to races in other moduli as well). In fact, the first value for which this occurs is at $608,981,813,029$ (discovered by Bays and Hudson in 1976). 

In 1962, Knapowski and Tur\'{a}n conjectured that if we randomly pick an arbitrarily large value of $x$, then there will ``almost certainly'' be more primes of the form $3n+2$ than $3n+1$ up to $x$. However, the Knapowski-Tur\'{a}n conjecture was later disproved by Kaczorowski and Sarnak, each working independently. In fact if we let $\nu$ denote the number of values of $x(\leq X)$ for which there are more primes of the form $3n+2$, the proportion $\frac{\nu}{X}$ does not tend to any limit as $X \rightarrow \infty$, but instead fluctuates. This opens the question of: what happens if we go out far enough? Will the race be unbiased if we set $X$ sufficiently far away from $0$? That is, is Chebyshev's bias only apparent for ``small'' values of X?

In 1994, while working with the mod $4$ race, Rubinstein and Sarnak introduced the logarithmic measure to find the percentage of time a certain team is in the lead \cite{chebyshevs bias}. Instead of counting 1 for each $x(\leq X)$ where there are more primes of the form $4n+3$ than of the form $4n+1$,  Rubinstein and Sarnak count $\frac{1}{x}$ . Instead of $\nu$, the sum is then approximately $\ln X$. They then scale this with the exact value of $\ln X$ to find the approximate proportion of time the $4n+3$ team is in the lead:  
\begin{equation*}
1 = \frac{\ln X}{\ln X} >
\left(\frac{1}{\ln X} \cdot \sum_{x \leq X} \frac{1}{x}\right) \rightarrow 0.9959 \ldots
\end{equation*}
where $x$ in the summation is only over values where there are more primes of the form $4n+3$ than of the form $4n+1$.

For the mod $3$ race, we have:
\begin{equation*}
\left(\frac{1}{\ln X} \cdot \sum_{x \leq X} \frac{1}{x}\right) \rightarrow 0.9990 \ldots
\end{equation*}
Using the logarithmic measure, we see that the $3n+2$ team is in the lead 99.9\% of the time!

\subsection{The Gaussian Primes}

\begin{defn}
A \textit{Gaussian integer} is a complex number whose real and imaginary parts are both integers. The Gaussian integers form an integral domain, which we denote with $\mathbb{Z}[i]$. 
\end{defn}

\begin{figure}[h]
\begin{center}
\includegraphics[width=9cm]{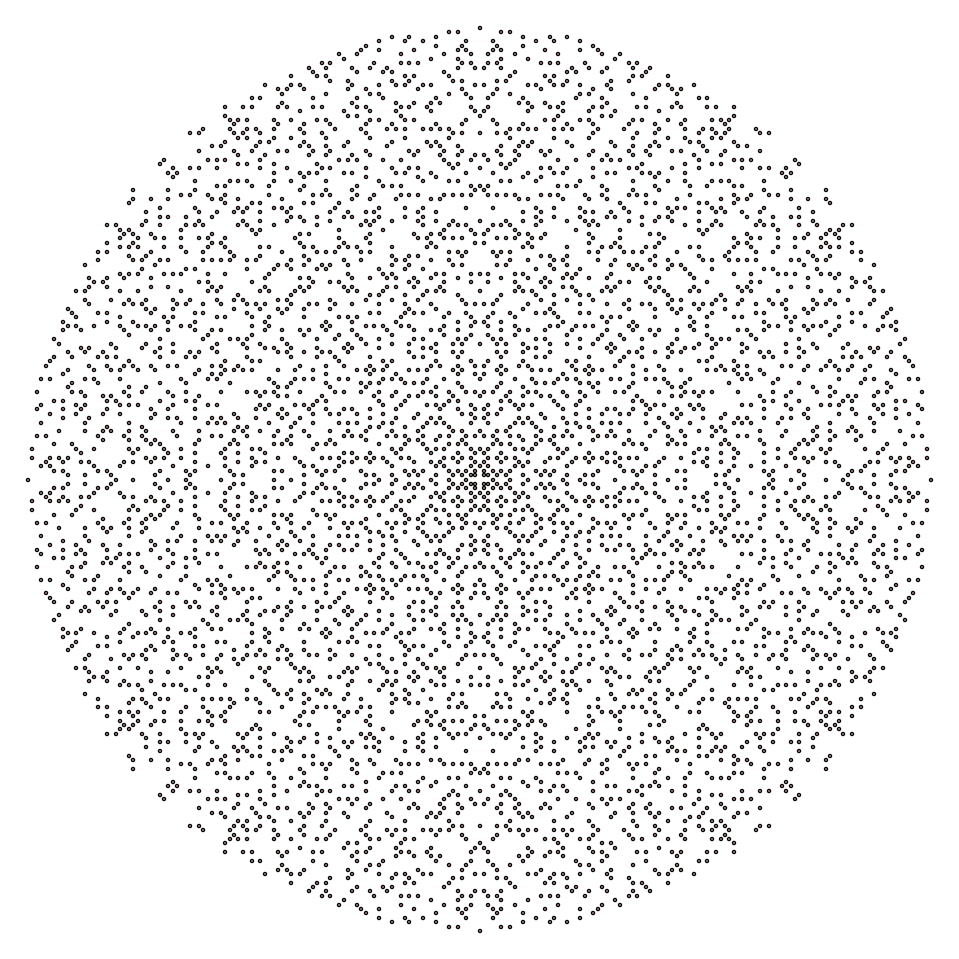}
\end{center}
\caption{Plot of Gaussian primes with norm $\leq 103^2$}
\label{gaussianprimes}
\end{figure}

In other words, for $i^2 = -1$, we have:
\begin{equation*}\mathbb{Z}[i] = \{a+bi | a,b \in \mathbb{Z}\}.\end{equation*}
The units of $\mathbb{Z}[i]$ are $\pm i$ and $\pm 1$. In addition, we say that two elements, $\mu$ and $\nu$ are \textit{associated} if $\mu=u\nu$ for $u$ being a unit in $\mathbb{Z}[i]$. Because of the four units, Gaussian primes (along with their complex conjugates) have an eightfold symmetry in the complex plane (figure \ref{gaussianprimes}). For convenience, we often write ``primes'' in place of ``primes unique up to associated elements.''

\begin{defn}
We say that an element in $\mathbb{Z}[i]$ is a \textit{Gaussian prime} if it is irreducible, i.e. if its only divisors are itself and a unit in $\mathbb{Z}[i]$.
\end{defn}

One might initially believe that the primes in $\mathbb{Z}$ are also irreducible elements in $\mathbb{Z}[i]$. However, this is not the case. In fact, there is a surprising connection between primes in mod $4$ arithmetic progressions in $\mathbb{Z}$ and the Gaussian primes. To understand this connection, we must first introduce the concept of norm. 

\begin{defn}
The \textit{norm} function takes a Gaussian integer $a+bi$ and maps it to a strictly positive real value. We denote the norm of a Gaussian integer as $N(a+bi) = (a+bi)(\overline{a+bi}) = (a+bi)(a-bi) = a^2 + b^2$. In other words, the norm function takes a Gaussian integer and multiplies it by its complex conjugate. One can geometrically understand the norm as the squared distance from the origin. 
\end{defn}

Let $\gamma = \alpha \cdot \beta$. The norm function is multiplicative; i.e. for $\gamma, \alpha, \beta$ elements in $\mathbb{Z}[i]$, 
\begin{equation*}N(\gamma) = N(\alpha \beta) = N(\alpha) N(\beta)\end{equation*}
We also note that the norm of any unit is 1. For example, if $\alpha = i = 0 + 1i$, then $N(\alpha) = 0^2 + 1^2 = 1$. In addition, we note that if an integer can be written as a sum of two squares, we can reduce it to two elements with smaller norms. For example, we note that $5=2^2 + 1^2 = (2+i)\cdot (2-i)= (2+i)\cdot (\overline{2+i})  = N(2+i)$. Thus, if a prime $p$ (in $\mathbb{Z}$) can be written as a sum of squares, we know it is not a prime element in $\mathbb{Z}[i]$. 

\begin{prpn}
\label{prpn: odd primes congruent to 1 mod 4}
If an odd prime $p$ is a sum of squares, it is congruent to $1 \pmod 4$ and not a prime element in $\mathbb{Z}[i]$.
\end{prpn}
Suppose $p = a^2 + b^2$. Since $p$ is odd, exactly one of $a$ or $b$ must be odd, and the other even. For the proof, we let $a$ be odd. Let $a = 2m+1$ and let $b=2n$. Then we have:
\begin{align*}
p & = a^2 + b^2 \\
&= (2m+1)^2 + (2n)^2 \\
&= 4m^2 + 4m +1 + 4n^2 \\
p & \equiv 1 \pmod 4
\end{align*}
Thus if $p \equiv 1 \pmod 4$, $p$ represents the \textit{norm} of two primes in $\mathbb{Z}[i]$. For example, $p = 13 \equiv 1 \pmod 4$ and $13 = N(\pi_1) = N(\pi_2)$, where  $\pi_1 = 2 +3i$ and $\pi_2 = 3 +2i$. We note that $\pi_2 = i \cdot \overline{\pi_1}$. (Here, we also note that counting primes in one quadrant is the same as counting primes unique up to associated elements).

\begin{prpn}
If an odd prime $p$ is congruent to $3 \pmod 4$, then $p$ is a prime element in $\mathbb{Z}[i]$.
\end{prpn}
For the proof, suppose for contradiction that we can factor $p$ into $(a+bi)\cdot(c+di)$. Using the multiplicative property of the norm function, we have:
\begin{align*}
N(p) &= N(a+bi)\cdot N(c+di) \\
p^2 & = (a^2 +b^2) \cdot (c^2 + d^2)
\end{align*}
Since $p$ is prime, $p^2$ can only be either $1 \cdot p^2$ or $p \cdot p$. Since we do not want a unit as a factor, we let $(a^2 + b^2) = p$ and $(c^2 + d^2) = p$. However, by proposition \ref{prpn: odd primes congruent to 1 mod 4}, we know that a solution would imply that $p$ is a sum of squares; i.e. $p \equiv 1 \pmod 4$. Thus, $p \equiv 3 \pmod 4$ cannot be factorized; i.e. $p$ is a Gaussian prime.

We now have enough information to classify a Gaussian prime into one of three general cases. Let $u$ be a unit in $\mathbb{Z}[i]$. Then:
\begin{itemize}
    \item{\makebox[2cm]{$u(1+i)$\hfill} Since $p = 2 = N(1+i)$}
    \item{\makebox[2cm]{$u(a+bi)$\hfill} $a^2 + b^2 = p \equiv 1 \pmod 4$}
    \item{\makebox[2cm]{$u(p)$\hfill} $p \equiv 3 \pmod 4$}
\end{itemize}

Thus, we can see that primes in $\mathbb{Z}$ with quadratic residue (modulo $4$) are not primes in $\mathbb{Z}[i]$. Instead, they represent the norms of two separate Gaussian primes. We can use this to derive an equation for the exact count of Gaussian primes (unique up to associated elements) within a certain norm. Let $\pi_G(x)$ represent the count of Gaussian primes up to norm x, then:
\begin{equation*}
\pi_G(x) = 2\pi(x;4,1) + \pi(\sqrt{x};4,3) + 1
\end{equation*}

The extra count is to include the Gaussian prime at $1+i$, which has norm $2$. 

In addition, we can extend our prime number theorem in the rational integers \eqref{eq:PNT} to a prime number theorem in the Gaussian integers by a modification of Dirichlet's Theorem \eqref{eq:PNTforAPs}. Moreover, we note the infinitude of Gaussian primes by their intimate connection with Dirichlet's Theorem for primes in mod $4$ progressions.

\begin{equation}
\label{eq: dirichlets gaussian thm}
\pi_G(x) \approx \frac{2x}{\phi(4)\log(x)} + \frac{\sqrt{x}}{\phi(4)\log(\sqrt{x})}
\end{equation}
The first term represents the approximation of primes congruent to $1$ (mod $4$), which are the norms of two primes in $\mathbb{Z}[i]$. The second term represents the approximation of primes congruent to $3$ (mod $4$), which have a norm of $p^2$ for $p \in 4n+3$. More precisely, we have:
\begin{equation*}
\lim_{x \rightarrow \infty} \dfrac{\pi_G(x)}{\frac{2x}{\phi(4)\log(x)} + \frac{\sqrt{x}}{\phi(4)\log(\sqrt{x})}} = 1
\end{equation*}

\pagebreak

The following code can be used in Sage to generate plots of Gaussian primes within a specified norm\footnote{We also created a video animation of Gaussian prime plots with norms from $10^1$ to $10^7$: \href{https://youtu.be/jRBCmXGlVJU}{https://youtu.be/jRBCmXGlVJU}}:

\begin{verbatim}
def gi_of_norm(max_norm):
    Gaussian_primes = {}
    Gaussian_integers = {}
    Gaussian_integers[0] = [(0,0)]
    for x in range(1, ceil(sqrt(max_norm))):
        for y in range(0, ceil(sqrt(max_norm - x^2))):
            N = x^2 + y^2
            if Gaussian_integers.has_key(N):
                Gaussian_integers[N].append((x,y))
            else:
                Gaussian_integers[N] = [(x,y)]
            if(y == 0 and is_prime(x) and x%4==3):
                have_prime = True
            elif is_prime(N) and N%4==1 or N==2:
                have_prime = True
            else:
                have_prime =False
            if have_prime:
                if Gaussian_primes.has_key(N):
                    Gaussian_primes[N].append((x,y))
                else:
                    Gaussian_primes[N] = [(x,y)]
    return Gaussian_primes,Gaussian_integers
def all_gaussian_primes_up_to_norm(N):
    gips = gi_of_norm(N)[0]
    return flatten([uniq([(x,y), (-y,x), (y,-x), (-x,-y)]) for x,y in flatten(gips.values(),
    max_level=1)], max_level=1)
N=10609 + 1 ### Declare norm here (in place of 10609)
P=scatter_plot(all_gaussian_primes_up_to_norm(N), markersize=RR(1000)/(N/50))
P.show(aspect_ratio=1, figsize=13, svg=False, axes = False)
\end{verbatim}
\section{Findings in the Rational Primes}
\subsection{Bias in the Legendre Symbols of Primes Modulo Another Prime}
\label{Legendre Symbol Races}
One phenomenon we wished to study in detail was Chebyshev's bias, specifically in regards to a randomly selected prime being more likely to have nonquadratic residue modulo some other prime. We approached this by first attempting to model the bias as a ``random'' walk using Legendre symbol values as steps. 

Let $q$ and $p$ be two randomly selected prime numbers. Then, according to Chebyshev's bias, $\left(\frac{q}{p}\right)$ has a slightly less than half probability of being a quadratic residue (i.e. returning a $1$). If we fix $p$ and let $q$ iterate through all primes, we get a sequence of $1$s and $(-1)$s (with the exception of when $q=p$, in which case we have $0$). If modeling as a random walk, the summation of our sequence should not wander far from $y=\sqrt{t}$, where $t$ denotes the index of the prime number $q$. Indeed, this is the case with all observed values of $p$ up to the final value of $q$ (we tested for primes $p < 1000$ and for $q$ iterating over primes $<10,000,000$). However, there is a noticeable bias in the summation. Most of the time, the summation of Legendre symbol values is negative, supporting the claim that there are slightly more nonquadratic residues. 
\pagebreak
\begin{figure}[h]
\begin{center}
\includegraphics[width=10cm]{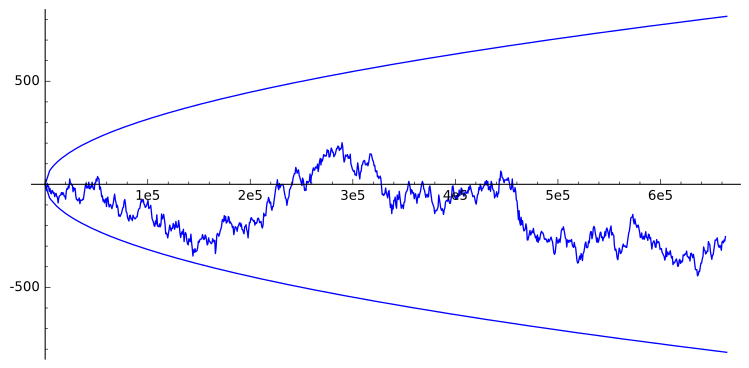}
\end{center}
\caption{Legendre Symbol walk for $p=97$}
\label{walk97all}
\end{figure}

We wished to model the average behavior of our Legendre symbol walks. To do this, we recorded the ratio of quadratic residues in each of our walks for $p$ fixed as we increase the range of primes over which $q$ iterates. For example, when $p$ = 97 and $q$ iterates over all primes less than $1000$, the ratio is $0.4698795$. When we allow $q$ to iterate over all primes less than $10,000,000$, the ratio of quadratic residues increases to $0.4997826$. We then plotted the average ratio for $167$ values of $p$ ($p \in \{3 \leq$ all primes $<1000\}$). In addition, we plotted the within-$p$ standard deviation of our ratio for each range of $q$ iterated. Since most primes have nonquadratic residue modulo another prime, the average ratio seems to converge to 0.50 from below as we increase the $q$-range.

\begin{figure}[h]
\begin{center}
\includegraphics[width=10cm]{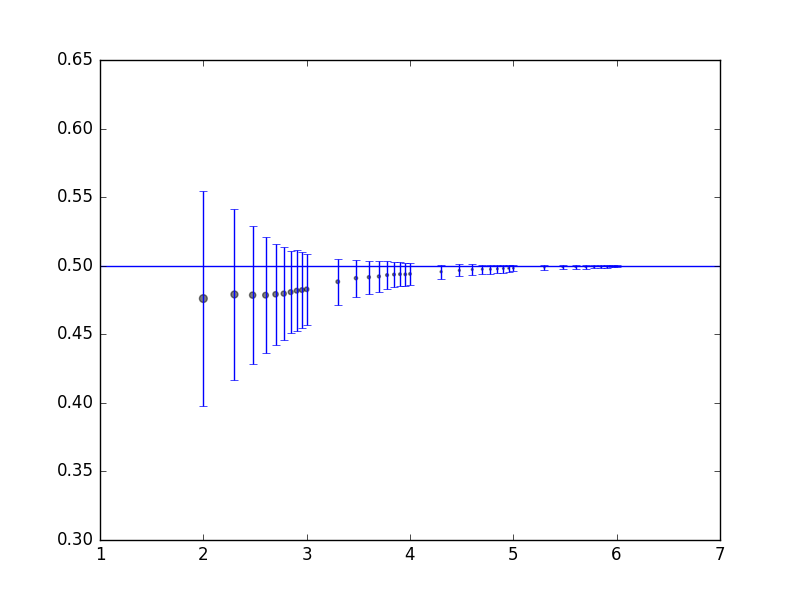}
\end{center}
\captionsetup{justification=centering,margin=2cm}
\caption{Plot of the Average Ratios. Horizontal axis denotes $\log(x)$, where $x$ is the range over which $q$ iterates. Vertical bars represent $1$ standard deviation}
\label{AvgRatioAll}
\end{figure}

We repeated our experiments with $\left(\frac{p}{q}\right)$ for $p$ fixed and $q$ varying and arrived at similar results. For $p \equiv 1 \pmod 4$, we know from quadratic reciprocity that $\left(\frac{p}{q}\right) = \left(\frac{q}{p}\right)$, so the contribution is the same (see theorem 10 in \cite{quadratic reciprocity clark}). For $p \equiv 3 \pmod 4$, $\left(\frac{p}{q}\right) \not = \left(\frac{q}{p}\right)$. However, Chebyshev's bias still exists (i.e. there are slightly fewer +1s than -1s). As a result, the average behavior is similar.

\subsection{Bias in the Legendre Symbols of consecutive Primes}
\label{bias in consecutive primes}

Our next experiment in the rational primes was to examine the ratio of consecutive quadratic or nonquadratic residues for primes $q$ modulo a fixed prime $p$. I.e. we wished to model the behavior of the ratio of $1,1$s or $-1, -1$s.

Since the probability of $q \pmod p$ being is quadratic residue is very slightly less than 0.5, we should expect expect our average ratio to converge to $\left(\frac{1}{2}\right)^{n-1}$ from below, where $n$ denotes the length of the consecutive chain. For example, for the ratio of three consecutive quadratic or nonquadratic residues, we expect to obtain approximately: $(\frac{1}{2})^3 + (\frac{1}{2})^3 \approx (\frac{1}{2})^{3-1}$. (The first term in the summation represents the probability of $3$ consecutive quadratic residues, and the second term represents the probability of $3$ consecutive nonquadratic residues). However, in a very recent paper (March, 2016), R. Lemke Oliver and K. Soundararajan \cite{unexpected_bias}, note that there is a much stronger bias in the residue of consecutive primes than expected. We set out to model this (stronger) bias with our Legendre symbol walk.

We repeated our average ratio experiment as in section \ref{Legendre Symbol Races}. However, we instead searched for $2$, $3$, and $4$ consecutive residues having the same sign. We notice that the average ratios converge to their expected values quite slowly, supporting R. Lemke Oliver and K. Soundararajan's recent discovery.

\begin{figure}[htp]
\centering
\includegraphics[width=.33\textwidth]{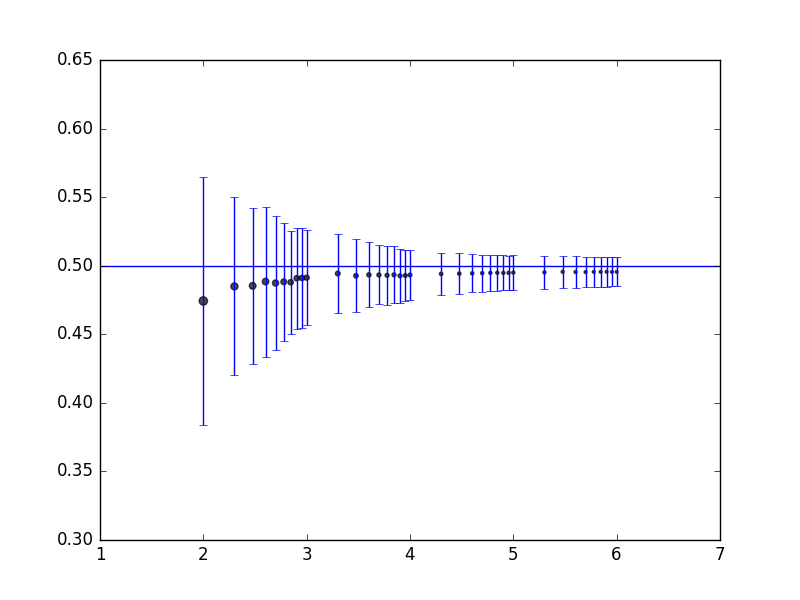}\hfill
\includegraphics[width=.33\textwidth]{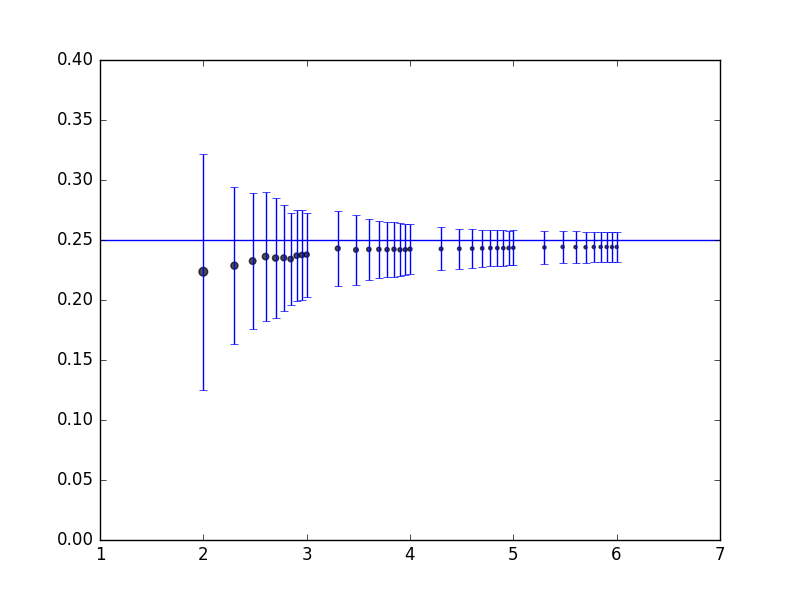}\hfill
\includegraphics[width=.33\textwidth]{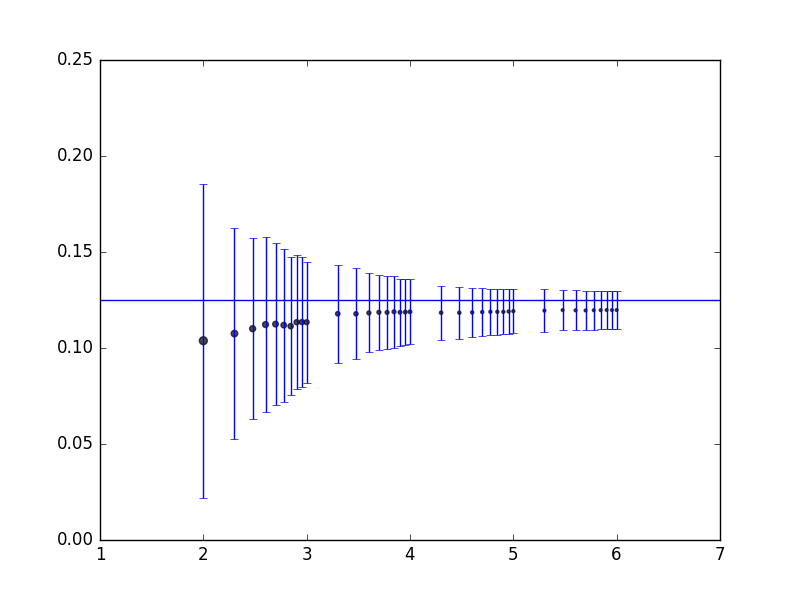}
\caption{From Left to Right: Two Consecutive, Three Consecutive, Four Consecutive}
\label{fig:AvgConsecutives}
\end{figure}

\subsection{Bias in the Legendre Symbols Modulo Primes in the Mod 4 Races}
\label{mod4 races}

We repeated our Legendre symbol walk with fixed $p$, but for $q$ varying only over primes congruent to $1$ (mod $4$), and again with primes congruent to $3$ (mod $4$). We observed Chebyshev's bias in both cases (on average). However, when $q$ varied over primes congruent to $1$ (mod $4$), we noticed a much stronger bias. For example, if we consider the walks for $p$=97, the walk for $q \equiv 1 \pmod 4$ seems to lie mostly below the $t$-axis. On the other hand, the walk for $q \equiv 3 \pmod 4$ seems to lie mostly above the $t$-axis. 

\begin{figure}[htp]
\centering
\includegraphics[width=.33\textwidth]{pQ97_p1000000_all.png}\hfill
\includegraphics[width=.33\textwidth]{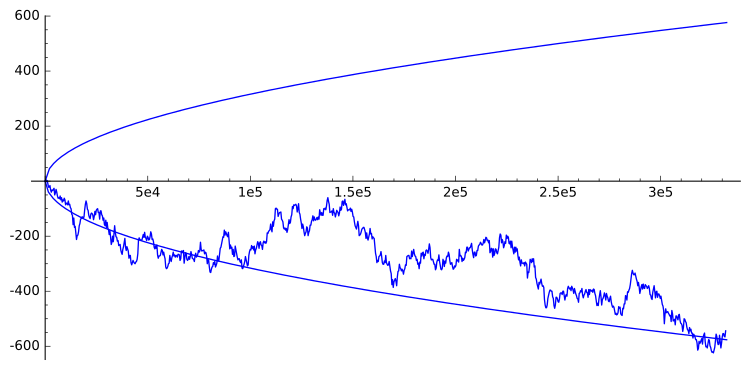}
\includegraphics[width=.33\textwidth]{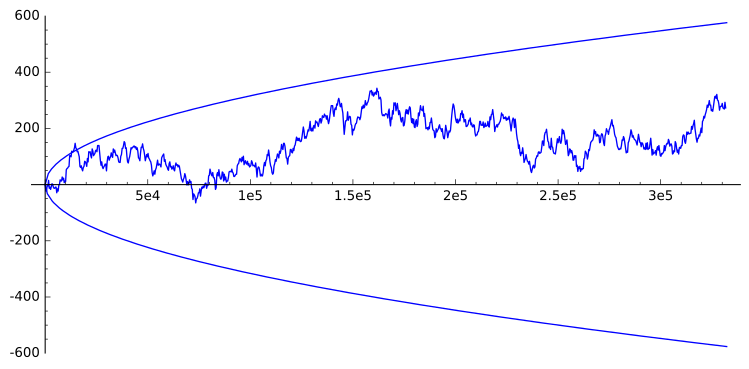}\hfill
\captionsetup{justification=centering,margin=2cm}
\caption{From Left to Right: iterating over all $q$, $q \equiv 1$ (mod $4$), $q \equiv 3$ (mod $4$). Iteration range for $q$ in all plots is $10,000,000$}.
\label{fig:97walks}
\end{figure}

We wished to check if this pattern exists on average. For $q\equiv 1$ (mod $4$), the average converges to $0.50$ more slowly than the average for $q\equiv 3$ (mod $4$). It seems that only allowing $q$ to iterate over primes with nonquadratic residue (mod $4$) removes, or at least diminishes, some part of Chebyshev's bias. We noticed a similar, but less distinct (see section \ref{bias in consecutive primes} and  [\href{unexpected_bias}{7}]), pattern while testing for consecutive residues being the same.

\begin{figure}[htp]
\centering
\includegraphics[width=.5\textwidth]{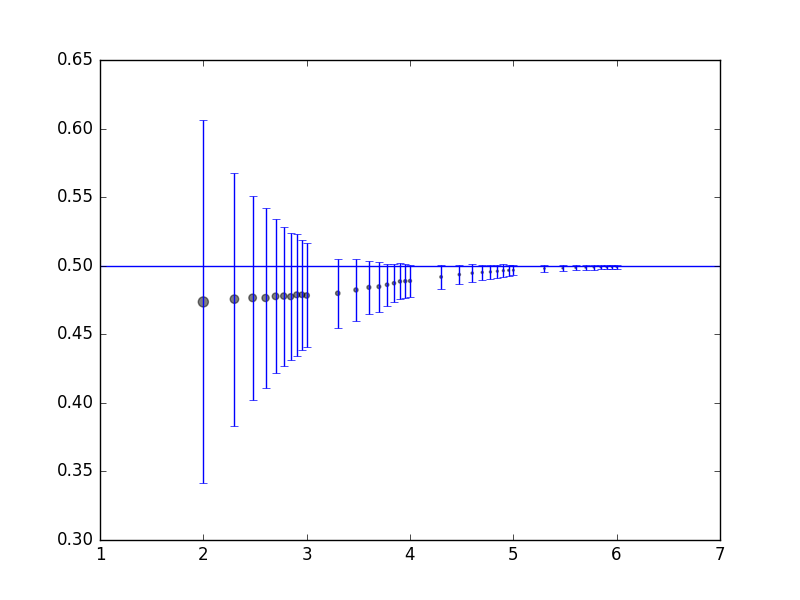}\hfill
\includegraphics[width=.5\textwidth]{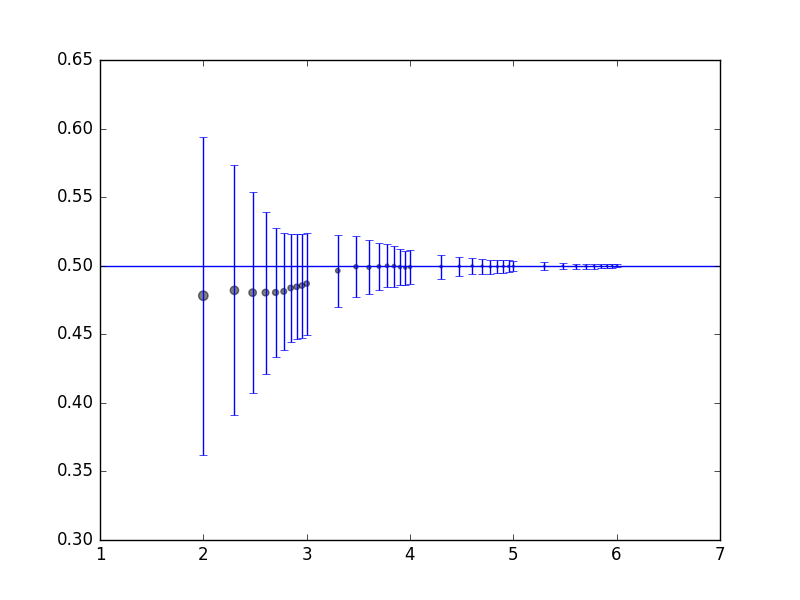}
\captionsetup{justification=centering,margin=2cm}
\caption{Average ratios of quadratic residues for $\left(\frac{q}{p}\right)$ \\ Left: $q \equiv 1$ (mod $4$). \\ Right:$q \equiv 3$ (mod $4$)}.
\label{fig:QuadNonquadratios}
\end{figure}

The following simple code can be used in Sage to generate a plot for Legendre symbol walks of $\left(\frac{q}{p}\right)$:

\begin{verbatim}
#declares maximum q-iteration range
maxN=10^7 
#P must be an odd prime for legendre_symbol(q,P) to be defined
P = 97 

primes = prime_range(3, maxN)
pm4={1:[], 3:[]}
pm4[1] = [q for q in primes if q % 4 == 1]
pm4[3] = [q for q in primes if q % 4 == 3]

#replace "3" with "1" to model walk with quadratic residues (mod 4)
lqP = [legendre_symbol(q, P) for q in pm4[3]] 

print "Legendre symbol walk for P={} and q iterating over primes less than {}".format(P,maxN)

sum_lqP = TimeSeries(lqP).sums()
#replace "3" with "1" to model walk with quadratic residues (mod 4)
sum_lqP.plot()+plot([sqrt(x),-sqrt(x)],(x,0,len(pm4[3]))) 

\end{verbatim}
\pagebreak

\section{Findings in the Gaussian Primes}

Chebyshev's bias in the rational primes has been well-documented. However, there has been comparatively less experimental research on such a bias in the Gaussian primes. In this section, we extend our model of Legendre symbol walks to the Gaussian primes to see if a similar bias occurs. To do this, we must first introduce a way to map a Gaussian integer to its residue in the rational integers modulo a Gaussian prime.
\begin{prpn}
\label{prpn: isomorphism}
A map that sends a Gaussian prime $a+bi$ to a residue $r$ (mod $\pi$), where $\pi = \alpha + \beta i$, is an isomorphism of rings between $\mathbb{Z}[i]/\pi\mathbb{Z}[i]$ and $\mathbb{Z}/p\mathbb{Z}$, where $p = N(\pi)$. In particular, if $\pi$ is an irreducible element in $\mathbb{Z}[i]$, then the residue class ring $\mathbb{Z}[i]/\pi\mathbb{Z}[i]$ is a finite field with $N(\pi)$ elements. 
\end{prpn}
A rigorous proof of proposition \ref{prpn: isomorphism} can be found in \cite{quadratic_reciprocity} as Theorem 12.

We first start with a ``soft'' proof as motivation for calculating a residue before showing a more rigorous proof. For two primes $p$ and $q$, the Euclidean algorithm shows that the $\gcd(p,q) = 1$. This fact allows us to easily calculate the residue of $q \pmod p$. Let $p$ and $q$ be prime numbers with $q>p$. Let n and r be integers:
\begin{align*}
q &= pn + r \\
q-r &= pn \\
q-r &\equiv 0 \pmod p \\
r &\equiv q \pmod p
\end{align*}
Where $r$ is a element from $(\mathbb{Z}/p\mathbb{Z})^\times$; i.e. $r$ is an element from the set of totatives of $p$.

We can extend this algorithm to the Gaussian primes. Let $a+bi$ and $\pi = \alpha + \beta i$ denote Gaussian primes with $N(a+bi) > N(\alpha+\beta i) = N(\pi)$. We can then write:
\begin{align*}
a+bi &= \pi(\phi+i\psi) + r \\
a+bi &= (\alpha+\beta i)(\phi+i\psi) + r \\
a+bi &= \alpha \phi + \alpha i \psi + \beta i \phi - \beta \psi + r
\end{align*}
We then group the real and imaginary terms:
\begin{align*}
a &= \alpha \phi - \beta \psi + r \\
b &= \alpha \psi + \beta \phi
\end{align*}
Use the imaginary component to solve for $\psi$, then solve for $a$:
\begin{align*}
\psi &= \frac{b - \beta \phi}{\alpha} \\
a &= \alpha \phi - \beta\left(\frac{b-\beta \phi}{\alpha}\right) + r\\
a &= \alpha \phi - \frac{b\beta}{\alpha} + \frac{\beta^2 \phi}{\alpha} + r
\end{align*}
Rearrange, multiply both sides by $\alpha$, and solve for $r$:
\begin{alignat*}{3}
a + \frac{b\beta}{\alpha} + r &= \alpha \phi + \frac{\beta^2 \phi}{\alpha} \\
a\alpha + b\beta - r\alpha &= \phi(\alpha^2 + \beta^2) \\
a\alpha + b\beta -r\alpha &\equiv 0 & \pmod{\alpha^2 + \beta^2} \\
a\alpha + b\beta &\equiv r\alpha & \pmod{\alpha^2 + \beta^2}\\
r &\equiv a + \alpha^{-1}b\beta & \pmod{\alpha^2 + \beta^2} \numberthis \label{eq: residue}
\end{alignat*}
where $r$ is an element from $(\mathbb{Z}/(\alpha^2 + \beta^2)\mathbb{Z})^\times = (\mathbb{Z}/N(\pi)\mathbb{Z})^\times = (\mathbb{Z}/p\mathbb{Z})^\times$ since $\alpha^2 + \beta^2 = N(\pi) = p$. 
\pagebreak

The idea is to use this residue to calculate the value of a  \textit{Gaussian Legendre symbol} $\left[\frac{a+bi}{\pi}\right]$  with the hope of observing a bias as in the rationals. First, we must lay the groundwork by introducing several concepts. (A comprehensive reference by Nancy Buck regarding Gaussian Legendre symbols, which includes the full proofs for the following propositions, can be found in \cite{quadratic_reciprocity}. Since many of the proofs are quite lengthy, we will only highlight sections relevant for our model).

\begin{defn}
For $k, l, \pi \in \mathbb{Z}[i]$, let $\pi$ be a Gaussian prime $\not = u(1+i)$ and such that $k$ and $l$ are not divisible by $\pi$. The \textit{Gaussian Legendre symbol} has the following properties
\end{defn}
\begin{itemize}
\item $\left[\dfrac{k}{\pi}\right] = \left[\dfrac{l}{\pi}\right]$ for $k \equiv l \pmod \pi$
\item $\left[\dfrac{k}{\pi}\right] \cdot \left[\dfrac{l}{\pi}\right] = \left[\dfrac{kl}{\pi}\right]$
\end{itemize}
For $p=N(\pi)$, the second point can be equivalently expressed as:
\begin{equation*}
k^{\frac{p-1}{2}}l^{\frac{p-1}{2}} = (kl)^{\frac{p-1}{2}} \equiv \left[\dfrac{kl}{\pi}\right] \pmod \pi
\end{equation*}
In addition, we have an analog of Euler's criterion in the Gaussian Legendre symbols:
\begin{equation*}
\left[\dfrac{k}{\pi}\right] \equiv k^{(p-1)/2}
\end{equation*}

\begin{thm}
\label{thm: GLS}
Every Gaussian Legendre symbol can be expressed in terms of a Legendre symbol in the rational integers.
\end{thm}
In particular, we have the following two equations for $\left[\dfrac{k}{\pi}\right]$. Let $k = a+bi$, $\pi = \alpha + \beta i$, and $N(\pi) = p$. Then:
 \begin{alignat}{3}
 \label{eq:no imaginary part}
 \left[\dfrac{a+bi}{\alpha}\right] &= \left(\dfrac{a^2 + b^2}{\alpha}\right); \quad && \pi \equiv 3 & \pmod 4\\
 \label{eq:with imaginary part}
 \left[\frac{a+bi}{\alpha+\beta i}\right] &= \left(\frac{a\alpha + b\beta}{p}\right); \quad && N(\pi) \equiv 1 & \pmod 4
 \end{alignat}
Recall that if $\pi$ is a prime element in $\mathbb{Z}[i]$, a zero imaginary part implies that $\pi = \alpha \equiv 3 \pmod 4$. For the proof of equation \eqref{eq:no imaginary part}, we must show that there exists an element $x \in \mathbb{Z}[i]$ such that $x^2 \equiv a + bi \pmod \alpha$ has a solution. We set $x = \phi + \psi i$ so that $\phi^2 - \psi^2 +2\phi\psi i \equiv a + bi \pmod \alpha$. Then we have the following two congruences by grouping real and imaginary terms:
\begin{alignat*}{3}
\phi^2 - \psi^2 &\equiv a \pmod \alpha \\
2\phi\psi &\equiv a \pmod \alpha
\end{alignat*}
We then square each congruence and add them together to get:
\begin{equation*}
\phi^4 + 2\phi^2\psi^2 + \psi^4 = (\phi^2 + \psi^2)^2 \equiv a^2+b^2 \pmod \alpha
\end{equation*}
It then suffices to check that there exists $\phi$ and $\psi \in \mathbb{Z}[i]$ such that both congruences have simultaneous solutions for the cases $a \not \equiv 0 \pmod \alpha$ and $a \equiv 0 \pmod \alpha$ (shown in \cite{quadratic_reciprocity}). Doing so shows that $\left[\dfrac{a+bi}{\alpha}\right]=1$ if and only if $\left(
\dfrac{a^2+b^2}{\alpha}\right) = 1$. In other words, we arrive at equation \eqref{eq:no imaginary part}: $\left[\dfrac{a+bi}{\alpha}\right]=\left(\dfrac{a^2+b^2}{\alpha}\right)$.
 
We now wish to consider the more interesting case when $N(\pi) \equiv 1 \pmod 4$; i.e. when $\pi = \alpha+\beta i$  for $\alpha,\beta \in \mathbb{Z}\backslash\{0\}$ and $\pi \not = (1+i)$. Let $\alpha$ be odd and $\beta$ be even. Let $k = a+bi$ with $a,b \in \mathbb{Z}$ and $\gcd(\pi, k)=1$. As above, we wish to determine if $x^2 \equiv a+bi \pmod \pi$ has a solution for $x \in \mathbb{Z}[i]$.

Recall that $p = N(\pi)$ is a prime congruent to $1 \pmod 4$. By proposition \ref{prpn: isomorphism}, we know the set of congruence class representatives modulo $\pi$ is $\{0, 1, 2, \ldots, p-1\}$. This allows us to only consider $x \in \mathbb{Z}$ when determining if $x^2 \equiv a+bi \pmod \pi$ has a solution.

We start by writing our congruence as an equivalence. The congruence $x^2 \equiv a+bi \pmod \pi$ is solvable if and only if there exists $x, \phi, \psi \in \mathbb{Z}$ such that:
\begin{align*}
x^2 - a - bi &= (\phi + \psi i )(\alpha+\beta i) \\
x^2 - a - bi &= \phi\alpha + \phi\beta i + \alpha \psi i - \beta \psi
\end{align*}
We then group the real and imaginary terms into separate equations:
\begin{align*}
x^2 - \alpha &= \phi \alpha - \beta \psi \\
-b & = \phi\beta + \alpha \psi 
\end{align*}
Then we multiply the real part by $\alpha$ and the imaginary part by $\beta$ and add:
\begin{align*}
x^2 - a\alpha &= \phi\alpha^2-\beta\psi\alpha \\
-b\beta &= \phi\beta^2 + \alpha\beta\psi \\
x^2\alpha - a\alpha-b\beta &=\phi\alpha^2+\phi\beta^2 \\
x^2\alpha - a\alpha-b\beta &=p\phi \\
x^2\alpha &= p\phi + a\alpha+b\beta
\end{align*}
Converting back to a congruence statement modulo $p$, we arrive at the following result:
\begin{equation*}
x^2\alpha \equiv a\alpha +b\beta \pmod p\\
\end{equation*}
\begin{equation}
\left(\dfrac{x^2 \alpha}{p}\right) = \left(\dfrac{a\alpha + b\beta}{p}\right) = \left(\dfrac{\alpha}{p}\right)\left(\dfrac{a+\alpha^{-1}b\beta}{p}\right) = \left(\dfrac{\alpha}{p}\right)\left(\dfrac{r}{p}\right)
\end{equation}
All that remains is to show that $\left(\dfrac{\alpha}{p}\right) = 1$. To do this, we use the law of quadratic reciprocity as described in proposition \ref{prpn: quadratic reciprocity}:
\begin{equation*}
\left(\dfrac{\alpha}{p}\right) = \mathlarger(-1)^{(\alpha-1)(p-1)/4}\left(\dfrac{p}{\alpha}\right)
\end{equation*}
Since $p \equiv 1 \pmod 4$, $p-1 \equiv 0 \pmod 4$. Thus, $\left(\dfrac{\alpha}{p}\right) = \left(\dfrac{p}{\alpha}\right)$. In addition, recall that $p = \alpha^2 + \beta^2$, so $p \equiv \beta^2 \pmod \alpha$. Thus, we can write $\left(\dfrac{\alpha}{p}\right) = \left(\dfrac{p}{\alpha}\right) = \left(\dfrac{\beta^2}{\alpha}\right)= \left(\dfrac{\beta}{\alpha}\right)\left(\dfrac{\beta}{\alpha}\right)$. It is then clear that regardless of the value of $\left(\dfrac{\beta}{\alpha}\right)$, we have  $\left(\dfrac{\alpha}{p}\right) =1$.

In conclusion, we arrive at equation \eqref{eq:with imaginary part}:
\begin{equation*}
\left[\dfrac{a+bi}{\alpha+\beta i}\right] = \left(\dfrac{a\alpha + b\beta}{p}\right) = \left(\dfrac{r}{p}\right)
\end{equation*}
\textbf{The Experiment.}

While implementing our random walk model on Sage, we decided to fix $\pi = \alpha+\beta i$ and let $a+bi$ iterate over Gaussian primes in the first quadrant sorted by increasing norm. In the case of $a+bi = a \equiv 3 \pmod 4$, the sorting is obvious. However, when $N(a+bi) = q \equiv 1 \pmod 4$, there are exactly two (distinct) Gaussian primes  with norm $q$ (we have $a+bi$ and $b+ai = i(\overline{a+bi})$, where $a^2 + b^2 = q$). When this is the case, we sort by the size of the real component. (For example, when $q= 17 = N(1+4i)$ and $N(4+i)$, we find the residue of $1+4i \pmod \pi$ first and then proceed to find the residue of $4+i \pmod \pi$). 

\pagebreak
When viewed individually, the resulting plots resemble the Legendre symbol walks in section \ref{Legendre Symbol Races}. However, we observe an interesting phenomenon when comparing walks that have the same $p = N(\pi_1) = N(\pi_2)$ where $\pi_1$ and $\pi_2$ are fixed with $a+bi$ iterating. We noticed for some $p$, the plots for $\pi_1$ and $\pi_2$ have strong positive correlation. For other $p$, the plots for $\pi_1$ and $\pi_2$ have strong negative correlation.

\begin{figure}[htp]
\centering
\includegraphics[width=.5\textwidth]{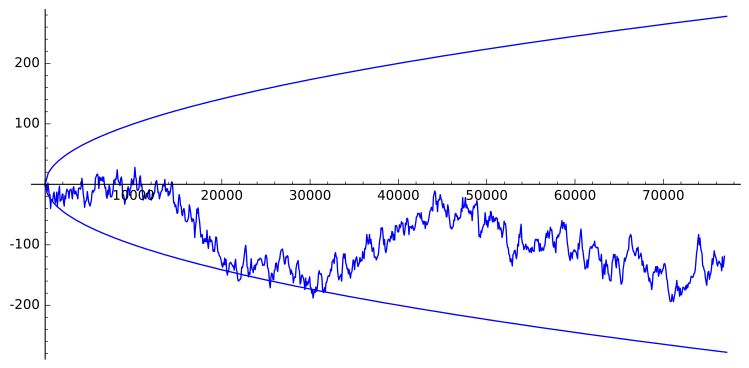}\hfill
\includegraphics[width=.5\textwidth]{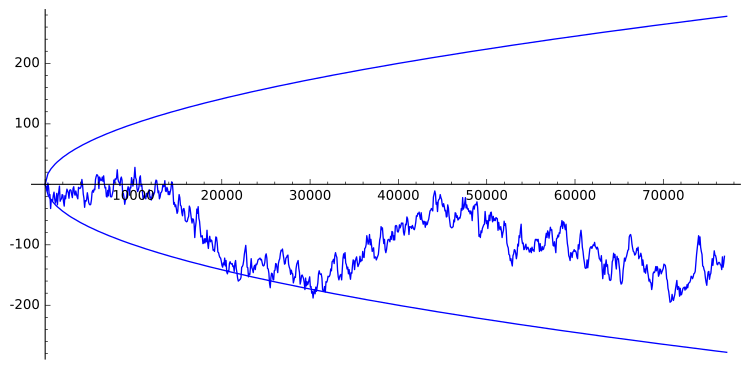}
\captionsetup{justification=centering,margin=2cm}
\caption{Gaussian Legendre symbol walks for $p=97$ \\ (strong positive correlation)\\ Left: $\left[\dfrac{a+bi}{4+9i}\right]$. \quad Right: $\left[\dfrac{a+bi}{9+4i}\right]$}.
\label{fig:97gaussianwalks}
\end{figure}

\begin{figure}[htp]
\centering
\includegraphics[width=.5\textwidth]{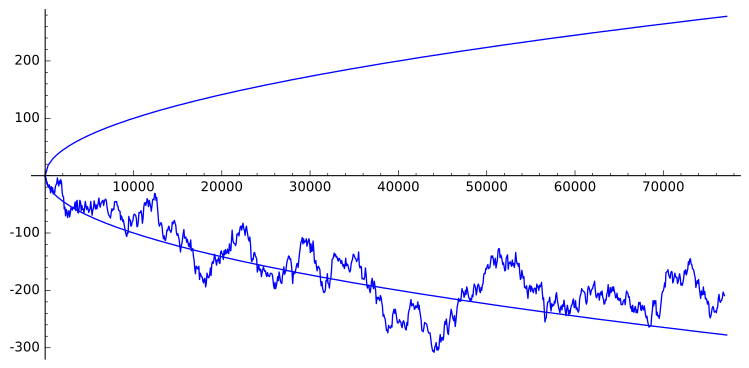}\hfill
\includegraphics[width=.5\textwidth]{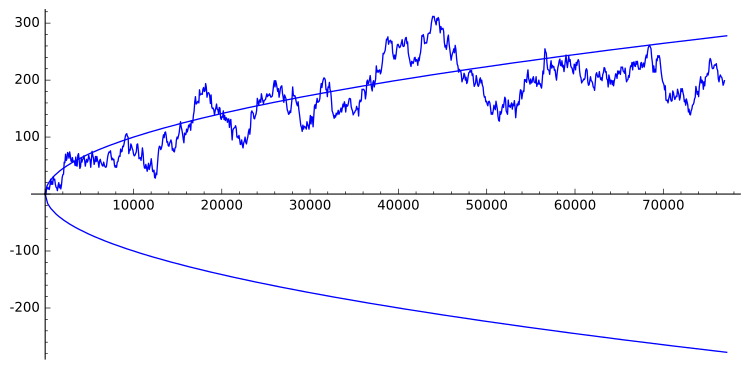}
\captionsetup{justification=centering,margin=2cm}
\caption{Gaussian Legendre symbol walks for $p=29$ \\ (strong negative correlation) \\ Left: $\left[\dfrac{a+bi}{2+5i}\right]$. \quad Right: $\left[\dfrac{a+bi}{5+2i}\right]$}.
\label{fig:29gaussianwalks}
\end{figure}

Before we attempt to (partially) explain this phenomenon, we must first introduce additional theory. 

\begin{thm}
The following $3$ properties hold for the Gaussian Legendre symbol
\end{thm}
\begin{align}
\label{eq: thm4 1}
\left[\dfrac{i}{\alpha+\beta i}\right] &= \left(-1\right)^{\frac{p-1}{4}} \\
\label{eq: thm4 2}
\left[\dfrac{1+i}{\alpha+\beta i}\right] &= \left(-1\right)^{\frac{(\alpha+\beta)^2-1}{8}} \\
\label{eq: thm4 3}
\left[\dfrac{a+bi}{\alpha+\beta i}\right] & =
\left[\dfrac{\alpha+\beta i}{a+bi}\right] 
\end{align}

\pagebreak
The proof of equation \eqref{eq: thm4 1} is as follows:

From Euler's criterion in the Gaussian Legendre symbols, we know that  $i^{(p-1)/2} \equiv \left[\dfrac{i}{\alpha+\beta i}\right] \pmod{\alpha+\beta i}$. We note that $i^{(p-1)/2}$ can be rewritten as follows:
\begin{equation*}
\mathlarger{i}^{\frac{p-1}{2}} = \mathlarger{i}^{2 \cdot \frac{p-1}{4}} = \mathlarger{(-1)}^{\frac{p-1}{4}}.
\end{equation*}
Thus, we have the congruence:
\begin{equation*}
\mathlarger{(-1)}^{\frac{p-1}{4}} \equiv \left[\dfrac{i}{\alpha+\beta i}\right]
\end{equation*}
For a proof by contradiction, we assume that the left side $\not \equiv$ the right side. Then let $-1 \equiv 1 \pmod{\alpha+\beta i}$. Converting the congruence to an equivalence, we get:
\begin{equation*}
-2 = (\alpha+\beta i)(\phi+\psi i)
\end{equation*}
We then take norms of both sides and simplify:
\begin{align*}
N(-2) &= N(\alpha + \beta i)N(\phi+\psi i) \\
4 &= p \cdot N(\phi + \psi i)
\end{align*}
This implies that $p \vert 4$, which cannot be true since $p \equiv 1 \pmod 4$. Therefore, we arrive at equation \eqref{eq: thm4 1}:
\begin{equation*}
\left[\dfrac{i}{\alpha+\beta i}\right] = \left(-1\right)^{\frac{p-1}{4}}.
\end{equation*}

For the proof of equation \eqref{eq: thm4 2}, we must consider two cases: when $\beta=0$ and when $\beta\not=0$.

Case 1: let $\beta=0$, so $p = \alpha^2$ and $\alpha \equiv 3 \pmod 4$. Recall our relations between the Gaussian Legendre symbols and the Legendre symbols in the rational integers as shown in theorem \ref{thm: GLS}. From equation \eqref{eq:no imaginary part}, we have:
\begin{equation*}
\left[\dfrac{1+i}{\alpha}\right] = \left(\dfrac{1+1}{\alpha}\right) = \left(\dfrac{2}{\alpha}\right)
\end{equation*}
Recall our second supplement of quadratic reciprocity in the rational integers. We can then express this as:
\begin{equation*}
\left(\dfrac{2}{\alpha}\right) = \mathlarger{(-1)}^{\frac{\alpha^2 -1}{8}} =  \mathlarger{(-1)}^{\frac{(\alpha+\beta)^2 -1}{8}}
\end{equation*}

Case 2: Let $\beta \not = 0$, so $p = \alpha^2 + \beta^2$ and $p \equiv 3 \pmod 4$. By equation \eqref{eq:with imaginary part}, we have:
\begin{equation*}
\left[\dfrac{1+i}{\alpha+\beta i}\right] = \left(\dfrac{\alpha+\beta}{p}\right).
\end{equation*}
Since our model only uses prime elements in the first quadrant, we assume that $\vert \alpha+\beta \vert > 1$ (the full proof without this assumption can be found in \cite{quadratic_reciprocity}). We continue by using the law of quadratic reciprocity:
\begin{equation*}
\left(\dfrac{\alpha+\beta}{p}\right) = \mathlarger{(-1)}^{\left(\frac{p-1}{2}\right)\left(\frac{\alpha+\beta-1}{2}\right)}\left(\dfrac{p}{\alpha+\beta}\right)
\end{equation*}
Since $p \equiv 1 \pmod 4$, then $\left(\frac{p-1}{2}\right)$ is always even. Thus, $\left(\frac{p-1}{2}\right)\left(\frac{\alpha+\beta -1}{2}\right)$ is even. So $\left(\dfrac{\alpha+\beta}{p}\right) = \left(\dfrac{p}{\alpha+\beta}\right)$.

\pagebreak
Next, we multiply $p$ by $2$ and apply a clever series of manipulations. We note that:
\begin{align*}
2p &= 2(\alpha^2 + \beta^2) \\
&=\alpha^2 + 2\alpha\beta + \beta^2 - 2\alpha\beta \\
&= (\alpha^2 + \beta^2)(\alpha^2-\beta^2) \\
0 &= (\alpha+\beta)^2 + (\alpha-\beta)^2 -2p \\
-(\alpha + \beta)^2 & = (\alpha-\beta)^2 -2p \\
0 & \equiv (\alpha-\beta)^2 - 2p \pmod{\alpha+\beta} \\
2p & \equiv (\alpha-\beta)^2 \pmod{\alpha+\beta}
\end{align*}
Let $x = (\alpha - \beta)^2$. Then there exists a solution to the congruence $x^2 \equiv 2p \pmod{\alpha+\beta}$. Then we have:
\begin{alignat*}{3}
\left(\dfrac{x^2}{\alpha+\beta}\right) &= \left(\dfrac{x}{\alpha+\beta}\right)\left(\dfrac{x}{\alpha+\beta}\right) &&= 1 \\
\left(\dfrac{x^2}{\alpha+\beta}\right) &= \left(\dfrac{2p}{\alpha+\beta}\right) &&= 1 \\
\left(\dfrac{2p}{\alpha+\beta}\right) &= \left(\dfrac{2}{\alpha+\beta}\right)\left(\dfrac{p}{\alpha+\beta}\right) &&= 1 
\end{alignat*}
Which implies that $\left(\dfrac{2}{\alpha+\beta}\right) = \left(\dfrac{p}{\alpha+\beta}\right) = \left(\dfrac{\alpha+\beta}{p}\right) = \left[\dfrac{1+i}{\alpha+\beta i}\right]$. Using the second supplement to quadratic reciprocity, we have:
\begin{equation*}
\left[\dfrac{1+i}{\alpha+\beta i}\right]  = \left(\dfrac{2}{\alpha+\beta}\right) =  \mathlarger{(-1)}^{\frac{(\alpha+\beta)^2-1}{8}}.
\end{equation*}
For the proof of equation \eqref{eq: thm4 3}, we must consider three cases:
\begin{enumerate}
\item $b  = \beta = 0$
\item $b = 0$ and $\beta \not = 0$ 
\item $b \not = 0$ and $\beta \not =0$.
\end{enumerate}
Case 1: Let $b = \beta = 0$. Then by equation \eqref{eq:no imaginary part}:
\begin{alignat*}{3}
\left[\dfrac{a}{\alpha}\right] &= \left[\dfrac{a^2}{\alpha}\right] &= 1\\
\left[\dfrac{\alpha}{a}\right] &= \left[\dfrac{\alpha^2}{a}\right] &= 1
\end{alignat*}
It is  then clear that $\left[\dfrac{a}{\alpha}\right] = \left[\dfrac{\alpha}{a}\right] = 1$.

Case 2: Assume $b = 0$ and $\beta \not = 0$. Then:
\begin{equation*}
\left[\dfrac{a}{\alpha+\beta i}\right] = \left(\dfrac{a\alpha}{p}\right) = \left(\dfrac{a}{p}\right)\left(\dfrac{\alpha}{p}\right) = \left(\dfrac{a}{p}\right)
\end{equation*}
(Recall we have already shown in theorem \ref{thm: GLS} that $\left(\dfrac{\alpha}{p}\right)=1$. Then we have:
\begin{equation*}
\left[\dfrac{\alpha+\beta i}{a}\right] = \left(\dfrac{\alpha^2+\beta^2}{a}\right) = \left(\dfrac{p}{a}\right)
\end{equation*}
From quadratic reciprocity, we know that $\left(\dfrac{a}{p}\right) = \mathlarger{(-1)}^{\frac{(p-1)(a-1)}{4}} \left(\dfrac{p}{a}\right)$
Since $p\equiv 1 \pmod 4$, we then see that $\left(\dfrac{a}{p}\right) = \left(\dfrac{p}{a}\right)$. Thus, we have:
\begin{equation*}
\left[\dfrac{a}{\alpha+\beta i}\right] = \left[\dfrac{\alpha+\beta i}{a}\right].
\end{equation*}

Case 3: Assume both $b$ and $\beta$ are nonzero. Since $a+bi$ and $\alpha+\beta i$ are distinct odd Gaussian primes, we have:
\begin{align*}
\left[\dfrac{a+bi}{\alpha+\beta i}\right] &= \left[\dfrac{a\alpha+b\beta}{p}\right] \\
\left[\dfrac{\alpha+\beta i}{a+ bi}\right] &= \left[\dfrac{a\alpha+b\beta}{q}\right]
\end{align*}
where $p = \alpha^2 + \beta^2$ and $q=a^2 + b^2$. Since we are working in the first quadrant, we assume that $a\alpha+b\beta >1$. We then wish to perform another manipulation (the idea is similar to the proof of equation \eqref{eq: thm4 2}). In particular, we wish to show that a certain congruence is solvable (mod $a\alpha+b\beta$). We note that:
\begin{alignat*}{3}
(a\alpha+b\beta)^2 + (a\beta-b\alpha)^2 & = a^2\alpha^2 +2ab\alpha\beta+b^2\beta^2+a^2\beta^2-2ab\alpha\beta+b^2\alpha^2 \\
&= a^2\alpha^2 +b^2\beta^2 +a^2\beta^2+ b^2\alpha^2 \\
&=(\alpha^2+\beta^2)(a^2+b^2) \\
(a\alpha+b\beta)^2 + (a\beta-b\alpha)^2 &= pq \\
(a\alpha+b\beta)^2 &= pq -(a\beta-b\alpha)^2 \\
0 &\equiv pq -(a\beta-b\alpha)^2 \pmod{a\alpha+b\beta}\\
pq &\equiv (a\beta-b\alpha)^2  \ \qquad \pmod{a\alpha+b\beta}
\end{alignat*}
We then set $a\beta-b\alpha = x$. Thus we have the congruence:
\begin{equation*}
pq \equiv x^2  \pmod{a\alpha+b\beta}.
\end{equation*}
To finish the proof, we show:
\begin{alignat*}{3}
\left(\dfrac{x^2}{a\alpha+b\beta}\right) &= \left(\dfrac{x}{a\alpha+b\beta}\right)\left(\dfrac{x}{a\alpha+b\beta}\right) &= 1 \\
= \left(\dfrac{pq}{a\alpha+b\beta}\right) &= \left(\dfrac{p}{a\alpha+b\beta}\right)\left(\dfrac{q}{a\alpha+b\beta}\right) &= 1
\end{alignat*}
which implies that $\left(\dfrac{p}{a\alpha+b\beta}\right) = \left(\dfrac{q}{a\alpha+b\beta}\right)$. Since we know that $p$ and $q$ are primes in $\mathbb{Z}$ that are congruent to $1 \pmod 4$, by quadratic reciprocity, we can equivalently write this as:  $\left(\dfrac{a\alpha+b\beta}{p}\right) = \left(\dfrac{a\alpha+b\beta}{q}\right)$. By applying equation \eqref{eq:with imaginary part} of theorem \ref{thm: GLS}, we then see that $\left[\dfrac{a+bi}{\alpha+\beta i}\right] = \left[\dfrac{\alpha+\beta i}{a+bi}\right]$.

\pagebreak

We now attempt to explain the strong ($\pm$) correlations we observed between Gaussian Legendre symbol walks with $\pi_1$ and $\pi_2$ fixed, where $\pi_2 = i\overline{\pi_1}$ and for $a+bi$ iterating over Gaussian primes in the first quadrant.

We first wish to establish a relationship between $\left[\dfrac{a+bi}{\alpha+\beta i}\right]$ and $\left[\dfrac{b+ai}{\alpha+\beta i}\right]$. This will allow us to find their combined contribution. (Recall the iteration order is one of $\left[\dfrac{a+bi}{\alpha+\beta i}\right] \rightarrow \left[\dfrac{b+ai}{\alpha+\beta i}\right]$ or $\left[\dfrac{b+ai}{\alpha+\beta i}\right] \rightarrow \left[\dfrac{a+bi}{\alpha+\beta i}\right]$, based on the size of the real part).

To find the conditions such that $\left[\dfrac{a+bi}{\alpha+\beta i}\right] = \left[\dfrac{b+ai}{\alpha+\beta i}\right]$ we set:
\begin{align*}
1 &= \left[\dfrac{a+bi}{\alpha+\beta i}\right]\cdot \left[\dfrac{b+ai}{\alpha+\beta i}\right] \\
&= \left[\dfrac{ab + a^2i +b^2i-ab}{\alpha+\beta i}\right] \\
&= \left[\dfrac{i}{\alpha+\beta i}\right]\cdot \left[\dfrac{a^2+b^2}{\alpha+\beta i}\right] \\
&=\mathlarger{(-1)}^{(p-1)/4}\left[\dfrac{q}{\alpha+\beta i}\right] \\
&= \mathlarger{(-1)}^{(p-1)/4}\left(\dfrac{q}{p}\right) \left(\dfrac{\alpha}{p}\right) \\
&= \mathlarger{(-1)}^{(p-1)/4}\left(\dfrac{q}{p}\right) \numberthis \label{eq: ab-ba relation}
\end{align*}
Thus, $\left[\dfrac{a+bi}{\alpha+\beta i}\right] = \left[\dfrac{b+ai}{\alpha+\beta i}\right]$ if $\frac{p-1}{4}$ is even and $\left(\dfrac{q}{p}\right) = 1$, or if $\frac{p-1}{4}$ is odd and $\left(\dfrac{q}{p}\right) = -1$. The conditions for the equivalence of $\left[\dfrac{a+bi}{\beta+\alpha i}\right] = \left[\dfrac{b+ai}{\beta+\alpha i}\right]$ are similar. 

Case 1: Let $\pi_1=\alpha+\beta i$ and $\pi_2 = \beta+ \alpha i$, where $N(\pi_1) = N(\pi_2) = p$. Let $\frac{p-1}{4}$ be an even integer. Suppose $\left(\dfrac{q}{p}\right) = 1$. Then by our equivalence relations, we have:
\begin{equation*}
\left[\dfrac{a+bi}{\alpha+\beta i}\right] = \left[\dfrac{b+ai}{\alpha+\beta i}\right] \text{and} \left[\dfrac{a+bi}{\beta+\alpha i}\right] = \left[\dfrac{b+bi}{\beta+\alpha i}\right]
\end{equation*}
Thus, whether the iteration order is $\left[\dfrac{a+bi}{\alpha+\beta i}\right] \rightarrow \left[\dfrac{b+ai}{\alpha+\beta i}\right]$ or $\left[\dfrac{b+ai}{\alpha+\beta i}\right] \rightarrow \left[\dfrac{a+bi}{\alpha+\beta i}\right]$, the combined contribution is one of $\pm2$. The same is true with $\left[\dfrac{a+bi}{\beta + \alpha i}\right] \rightarrow \left[\dfrac{b+ai}{\beta + \alpha i}\right]$ or $\left[\dfrac{b+ai}{\beta + \alpha i}\right] \rightarrow \left[\dfrac{a+bi}{\beta+\alpha i}\right]$.

We now consider the case when $\frac{p-1}{4}$ is still an even integer, but $\left(\frac{q}{p}\right) = -1$. Then by our equivalence relations, we have:
\begin{equation*}
\left[\dfrac{a+bi}{\alpha+\beta i}\right] \not = \left[\dfrac{b+ai}{\alpha+\beta i}\right] \text{and} \left[\dfrac{a+bi}{\beta + \alpha i}\right] \not = \left[\dfrac{b+ai}{\beta+\alpha i}\right]
\end{equation*}
Thus, for any norm-sorted iteration order, the  combined contribution will be 0. 

Case 2: Now we let $\frac{p-1}{4}$ be an odd integer. Suppose $\left(\dfrac{q}{p}\right) = 1$. From our equivalence relations, we know that:
\begin{equation*}
\left[\dfrac{a+bi}{\alpha+\beta i}\right] \not = \left[\dfrac{b+ai}{\alpha+\beta i}\right] \text{and} \left[\dfrac{a+bi}{\beta + \alpha i}\right] \not = \left[\dfrac{b+ai}{\beta+\alpha i}\right]
\end{equation*}
Then for any norm-sorted iteration order, the combined contribution from $a+bi$ and $b+ai$ will be zero for both walks of $\pi_2$ and $\pi_2$.

Now we consider the case when $\frac{p-1}{4}$ is still odd, but $\left(\frac{q}{p}\right) = -1$. In this case, $\left[\dfrac{a+bi}{\alpha+\beta i}\right] = \left[\dfrac{b+ai}{\alpha+\beta i}\right]$. Thus, for any norm-sorted iteration order, the combined contribution will be one of $\pm 2$.

If we can establish the conditions for equivalence between $\left[\dfrac{a+bi}{\alpha+\beta i}\right]$ and $\left[\dfrac{a+bi}{\beta+\alpha i}\right]$ we will be able to fully explain the strong positive and negative correlations observed. (Note: it still remains to show what happens when $a+bi$ iterates over Gaussian primes $a+bi = a\equiv 3 \pmod4$. However, since prime elements of this form are much more sparse by equation \eqref{eq: dirichlets gaussian thm}, we can ignore them for the purposes of our explanation). Unfortunately, we found it quite difficult to rigorously prove the equivalence conditions (in particular, because the Legendre (more precisely, Jacobi) symbol $\left(\dfrac{p}{\beta}\right)$ is not defined for $\beta$ an even integer), so we leave it as a conjecture.
\begin{conj}
The equivalence between $\left[\dfrac{a+bi}{\alpha+\beta i}\right]$ and $\left[\dfrac{a+bi}{\beta+\alpha i}\right]$ depends only on the value of the Legendre symbol $\left(\dfrac{q}{p}\right)$. In particular, $\left[\dfrac{a+bi}{\alpha+\beta i}\right] = \left[\dfrac{a+bi}{\beta+\alpha i}\right]$ if $\left(\dfrac{q}{p}\right) = 1$, and $\left[\dfrac{a+bi}{\alpha+\beta i}\right] \not= \left[\dfrac{a+bi}{\beta+\alpha i}\right]$ if $\left(\dfrac{q}{p}\right) \not= 1$.
\end{conj}

We will use the following shorthand notation for clarity and convenience:
\begin{align*}
\pi_{1a} &= \left[\dfrac{a+bi}{\alpha+\beta i}\right] \quad 
\pi_{1b} = \left[\dfrac{b+ai}{\alpha+\beta i}\right] \\
\pi_{2a} &= \left[\dfrac{a+bi}{\beta+\alpha i}\right] \quad
\pi_{2b} = \left[\dfrac{b+ai}{\beta+\alpha i}\right] \\
\pi_1 &= \pi_{1a} + \pi_{1b} \quad \quad \pi_2 = \pi_{2a} + \pi_{2b}
\end{align*}
To summarize, we have shown (conjectured) the following relations:
\begin{align}
\label{pi1a1b}
\pi_{1a} \pi_{1b} = \mathlarger{(-1)}^{(p-1)/4}\left(\frac{q}{p}\right) \\
\label{pi2a2b}
\pi_{2a} \pi_{2b} = \mathlarger{(-1)}^{(p-1)/4}\left(\frac{q}{p}\right) \\
\label{pi1a2a}
\pi_{1a} \pi_{2a} = \left(\frac{q}{p}\right) \\
\label{pi1b2b}
\pi_{1b} \pi_{2b} = \left(\frac{q}{p}\right)
\end{align}

We can now explain the strong ($\pm$) correlations between plots for $\pi_1$ and $\pi_2$ fixed. 

Consider the case when $\frac{p-1}{4}$ is even and $\left(\dfrac{q}{p}\right) = 1$. If $\pi_{1a} = 1$ (resp. $-1$), then by equation \eqref{pi1a1b}, $\pi_{1b}=1$ (resp. $-1$). Using equation \eqref{pi1a2a}, $\pi_{2a}=1$ (resp. $-1$), and by equation \eqref{pi2a2b}, $\pi_{2b} = 1$ (resp. $-1$). Thus, when $\frac{p-1}{4}$ is even and $\left(\dfrac{q}{p}\right) = 1$, the walks for $\pi_{1}$ and $\pi_{2}$ move exactly together with combined contribution one of $\pm 2$. Consider the case when $\frac{p-1}{4}$ is even and $\left(\dfrac{q}{p}\right) = -1$. If $\pi_{1a} = 1$ (resp. $-1$), then by equation \eqref{pi1a1b}, $\pi_{1b}=-1$ (resp. $1$). Using equation \eqref{pi1a2a}, $\pi_{2a}=-1$ (resp. $1$), and by equation \eqref{pi2a2b}, $\pi_{2b} = 1$ (resp. $-1$). Then $\pi_1$ and $\pi_2$ do not move together, but the combined contribution for that particular q is $0$, so there is little movement and the correlation remains close to $+1$.
\pagebreak

Consider the case when $\frac{p-1}{4}$ is odd and $\left(\dfrac{q}{p}\right) = 1$. If $\pi_{1a} = 1$ (resp. $-1$), then by equation \eqref{pi1a1b}, $\pi_{1b}=-1$ (resp. $1$). Using equation \eqref{pi1a2a}, $\pi_{2a}=1$ (resp. $-1$), and by equation \eqref{pi2a2b}, $\pi_{2b} = -1$ (resp. $1$). Thus, when $\frac{p-1}{4}$ is odd and $\left(\dfrac{q}{p}\right) = 1$, the walks move together, but with a combined contribution of $0$ for that particular $q$.  Consider the case when $\frac{p-1}{4}$ is odd and $\left(\dfrac{q}{p}\right) = -1$. If $\pi_{1a} = 1$ (resp. $-1$), then by equation \eqref{pi1a1b}, $\pi_{1b}=1$ (resp. $-1$). Using equation \eqref{pi1a2a}, $\pi_{2a}=-1$ (resp. $1$), and by equation \eqref{pi2a2b}, $\pi_{2b} = -1$ (resp. $1$). Then $\pi_1$ and $\pi_2$ move exactly \textit{opposite} to each other, causing the correlation to remain close to $-1$.

\section{Conclusions}

If one performs a Legendre symbol race in the rational primes, the sorting is obvious. However, if one extends the model to the Gaussian primes, the sorting is less clear. In this project, we only used one sorting order (by norm and then by size of real part). In addition, we only considered primes in the first quadrant. Perhaps future projects can model Gaussian Legendre symbol walks with different sorting orders, iterating over different combinations of quadrants, and up to greater norm values. Moreover, we mostly ignored the contribution of Gaussian primes of the form $a\equiv 3 \pmod4$ since they are much less numerous. Although it was not rigorously discussed, it seems that primes of this form contribute to a bias toward nonquadratic residues when comparing plots with odd $\frac{p-1}{4}$ (i.e. the plots with negative correlation). It would be interesting to quantify their effect on the correlation between the plots of $\pi_1$ and $\pi_2$. In addition, we noted in section \ref{mod4 races} that a Legendre symbol walk over rational primes $\equiv 3 \pmod 4$ seems to reduce some of Chebyshev's bias. It would be interesting to see an explanation for this phenomenon as well (perhaps there is an interesting connection to the Gaussian primes). We hope that we outlined enough theory for an inquisitive reader to begin asking their own questions about the fascinating Gaussian primes. 

\section*{Acknowledgments}

I would like to extend a special thank you to Dr. Stephan Ehlen for his guidance and teaching throughout the past year and through the duration of this project.

I would also like to thank  Dr. Henri Darmon of McGill University, le Centre de recherches mathématiques, and l'Institut des sciences mathématiques for providing me with funding and the opportunity to research this topic.

In addition, I would like to thank Dr. Yara Elias and Dr. Kenneth Ragan for their excellent teaching, and for helping me secure this research project. 

\pagebreak

\textsc{McGill University, Desautels Faculty of Management, 1001 Sherbrooke St. West, Montreal, Quebec, Canada H3A, 1G5} \\
\textit{E-mail address}, D. Hutama: \quad \href{mailto:daniel.hutama@mail.mcgill.ca}{daniel.hutama@mail.mcgill.ca}

\end{document}